\documentclass[a4paper]{article}

\usepackage{amssymb}

\def\GAP{\textsf{GAP}}
\def\ATLAS{\textsc{Atlas}}
\def\Co{\mathrm{Co}}
\def\E{\mathrm{E}}
\def\F{\mathrm{F}}
\def\Fi{\mathrm{Fi}}
\def\HN{\mathrm{HN}}
\def\He{\mathrm{He}}
\def\J{\mathrm{J}}
\def\McL{\mathrm{McL}}
\def\Suz{\mathrm{Suz}}
\def\Th{\mathrm{Th}}
\def\MAGMA{\textsf{MAGMA}}

\def\Sp{{\rm Sp}}
 \def\Z{{\mathbb Z}} \def\Q{{\mathbb Q}}

\def\B{{\mathbb B}}
\def\M{{\mathbb M}}
\def\tthdump#1{#1}
\tthdump{}

\begin{document}

\tthdump{\title{Some steps in the verification of the ordinary character table of the Monster group}}

\author{\textsc{Thomas Breuer, Kay Magaard, Robert A.~Wilson}}

\date{December 12th, 2024}

\maketitle

\abstract{%
We show the details of certain computations that are used in \cite{Mverify}.}

\textwidth16cm
\oddsidemargin0pt

\parskip 1ex plus 0.5ex minus 0.5ex
\parindent0pt

\tableofcontents


\section{Overview}

The aim of~\cite{Mverify} is to verify the ordinary character table
of the Monster group $\M$.
Here we collect,
in the form of an explicit and reproducible
{\GAP}~\cite{GAP} session protocol,
the relevant computations that are needed in that paper.

We proceed as follows.

Section~\ref{natural} verifies the decomposition of the restrictions of the
ordinary irreducible character of degree $196\,883$ of $\M$
to the subgroups $2.\B$ and $3.\Fi_{24}^{\prime}$ (and $3.\Fi_{24}$),
as stated in~\cite[Lemma~1]{Mverify}.

Section~\ref{suborbits} verifies the decompositions of
the transitive constituents of the permutation character
of the action of $C_{\M}(a) \cong 2.\B$ on the conjugacy class $a^{\M}$,
where $a$ is a \texttt{2A} involution in $\M$.

Sections~\ref{Mclasses} and~\ref{Mpowermaps} construct the
character table head of $\M$, that is,
the lists of conjugacy class lengths, element orders, and power maps.

Section~\ref{sect:natcharM} constructs the values of the irreducible
degree $196\,883$ character of $\M$
and decides the isomorphism type of the \texttt{3B} normalizer in $\M$.

With this information and with the (already verified) character tables
of the subgroups $2.\B$, $2^{1+24}_+.\Co_1$, $3.\Fi_{24}$, and
$3^{1+12}_+.2.\Suz.2$ of $\M$,
computing the irreducible characters of $\M$ is then easy;
this corresponds to~\cite[Section~5]{Mverify}
and is done in Section~\ref{sect:irreduciblesM}.

The final sections~\ref{sect:table_c2b},
\ref{sect:norm3B}, \ref{sect:table_N5B} document the constructions of three
character tables of subgroups of $\M$.

We will use the {\GAP} Character Table Library
and the interface to the {\ATLAS} of Group Representations~\cite{AGRv3},
thus we load these {\GAP} packages.

\begin{verbatim}
    gap> LoadPackage( "ctbllib", false );
    true
    gap> LoadPackage( "atlasrep", false );
    true
\end{verbatim}

The {\MAGMA} system~\cite{Magma} will be needed
for computing some character tables
and for many conjugacy tests.
If the following command returns \texttt{false}
then these steps will not work.

\begin{verbatim}
    gap> CTblLib.IsMagmaAvailable();
    true
\end{verbatim}

We set the line length to $72$, like in other standard testfiles.

\begin{verbatim}
    gap> SizeScreen( [ 72 ] );;
\end{verbatim}

\section{Some restrictions of the natural character of $\M$}%
\label{natural}

We assume the existence of an ordinary irreducible character $\chi$
of degree $196\,883$ of the Monster group $\M$,
and that $\M$ has only two conjugacy classes of involutions.

First we compute the restriction of $\chi$ to $2.\B$.

The only faithful degree $196\,883$ character of $2.\B$
that has at most two different values on involutions is
$1a + 4371a + 96255a + 96256a$, as claimed in \cite[Lemma~1]{Mverify}.
This follows from the following data about $2.\B$.

\begin{verbatim}
    gap> table2B:= CharacterTable( "2.B" );;
    gap> cand:= Filtered( Irr( table2B ), x -> x[1] <= 196883 );;
    gap> List( cand, x -> x[1] );
    [ 1, 4371, 96255, 96256 ]
    gap> inv:= Positions( OrdersClassRepresentatives( table2B ), 2 );
    [ 2, 3, 4, 5, 7 ]
    gap> PrintArray( List( cand, x -> x{ Concatenation( [ 1 ], inv ) } ) );
    [ [       1,       1,       1,       1,       1,       1 ],
      [    4371,    4371,    -493,     275,     275,      19 ],
      [   96255,   96255,    4863,    2047,    2047,     255 ],
      [   96256,  -96256,       0,    2048,   -2048,       0 ] ]
\end{verbatim}

Note that $96256a$ must occur as a constituent
because it is the only faithful candidate,
and it can occur only once because otherwise only $4371a + 2 \cdot 96256a$
or $4371 \cdot 1a + 2 \cdot 96256a$ would be possible decompositions,
which have more than two values on involution classes.
Thus the values of $\chi_{2.\B}$ on the classes 4 and 5 differ by $4096$.

If $96255a$ would {\bf not} occur then the values of $\chi_{2.\B}$
on the classes 3 and 7 would differ by $512$ times the multiplicity of $4371$,
but $65659 \cdot 1a + 8 \cdot 4371a + 96256a$ is not a solution.
Thus $96255a$ must occur exactly once.

The sum of $96255a$ and $96256a$ has four different values
on involutions, hence also $4371a$ must occur.

We see that the values of $\chi$ on the classes of involutions are
$4371$ and $275$, respectively.

\begin{verbatim}
    gap> Sum( cand ){ inv };
    [ 4371, 4371, 4371, 275, 275 ]
\end{verbatim}

The restriction of $\chi$ to $3.\Fi_{24}^{\prime}$ is computed similarly,
as follows.

Exactly seven irreducible characters of $3.\Fi_{24}^{\prime}$ can occur as
constituents of the restriction of $\chi$.

\begin{verbatim}
    gap> table3Fi24prime:= CharacterTable( "3.Fi24'" );;
    gap> cand:= Filtered( Irr( table3Fi24prime ), x -> x[1] <= 196883 );;
    gap> inv:= Positions( OrdersClassRepresentatives( table3Fi24prime ), 2 );
    [ 4, 7 ]
    gap> mat:= List( cand, x -> x{ Concatenation([1], inv)});;
    gap> PrintArray( mat );
    [ [      1,      1,      1 ],
      [   8671,    351,    -33 ],
      [  57477,   1157,    133 ],
      [    783,     79,     15 ],
      [    783,     79,     15 ],
      [  64584,   1352,     72 ],
      [  64584,   1352,     72 ] ]
\end{verbatim}

Since $\chi$ is rational, we need to consider only rationally irreducible
characters, that is, the possible constituents are $1a$, $8671a$, $57377a$,
$783ab$, and $64584ab$.

\begin{verbatim}
    gap> List( cand, x -> x[2] );
    [ 1, 8671, 57477, 783*E(3), 783*E(3)^2, 64584*E(3), 64584*E(3)^2 ]
\end{verbatim}

We see that the value on the first class of involutions must be $4371$,
since all values of the possible constituents are positive
and too large for the other possible value $275$.

Since the values of all possible constituents on the second class of
involutions are at most equal to the values on the first class,
and equal only for $1a$,
we conclude that the value on the second class of involutions is $275$.

We see from the ratio of the value on the identity element
and on the first class of involutions
that constituents of degree $57477$ or $2 \cdot 64584$ exist.

\begin{verbatim}
    gap> Float( 196883 / 4371 );
    45.043
    gap> List( mat, v -> Float( v[1] / v[2] ) );
    [ 1., 24.7037, 49.6776, 9.91139, 9.91139, 47.7692, 47.7692 ]
    gap> Float( ( 196883 - 2 * 64584 ) / ( 4371 - 2 * 1352 ) );
    40.6209
    gap> Float( ( 196883 - 57477 ) / ( 4371 - 1157 ) );
    43.3746
    gap> Float( ( 196883 - 2*57477 ) / ( 4371 - 2*1157 ) );
    39.8294
    gap> Float( ( 196883 - 3*57477 ) / ( 4371 - 3*1157 ) );
    27.1689
\end{verbatim}

First suppose that $64584ab$ is not a constituent.
The above ratios imply that (at least) three constituents of degree $57477$
must occur.

However, then the degree admits at most two constituents of degree $8671$,
hence the value on the second class of involutions cannot be $275$,
a contradiction.

This means that both $64584ab$ and $57477a$ occur with multiplicity one.

\begin{verbatim}
    gap> mat[3] + mat[6] + mat[7];
    [ 186645, 3861, 277 ]
\end{verbatim}

The second involution class forces one constituent of degree $8671$
(which is the only candidate that can contribute a negative value),
and then a character of degree $1567$ remains to be decomposed.
The only solution for the degrees of its constituents is $1 + 1566$.
We get the decomposition
$1a + 8671a + 57477a + 783ab + 64584ab$,
as claimed in \cite[Lemma~2]{Mverify}.

\begin{verbatim}
    gap> Sum( mat );
    [ 196883, 4371, 275 ]
\end{verbatim}

The characters of the degrees $1$, $8671$, and $57477$ extend two-fold
from $3.\Fi_{24}^{\prime}$ to $3.\Fi_{24}$.
In order to decompose the restriction of $\chi$ to $3.\Fi_{24}$,
we have to determine which extensions from $3.\Fi_{24}^{\prime}$ occur.
The following irreducible characters of $3.\Fi_{24}$ can occur as
constituents of the restriction of $\chi$.

\begin{verbatim}
    gap> table3Fi24:= CharacterTable( "3.Fi24" );;
    gap> cand:= Filtered( Irr( table3Fi24 ), x -> x[1] <= 196883 );;
    gap> inv:= Positions( OrdersClassRepresentatives( table3Fi24 ), 2 );
    [ 3, 5, 172, 173 ]
    gap> mat:= List( cand, x -> x{ Concatenation([1], inv)});;
    gap> PrintArray( mat );
    [ [       1,       1,       1,       1,       1 ],
      [       1,       1,       1,      -1,      -1 ],
      [    8671,     351,     -33,    1495,     -41 ],
      [    8671,     351,     -33,   -1495,      41 ],
      [   57477,    1157,     133,    5865,     233 ],
      [   57477,    1157,     133,   -5865,    -233 ],
      [    1566,     158,      30,       0,       0 ],
      [  129168,    2704,     144,       0,       0 ] ]
\end{verbatim}

We get the decomposition
$1a + 8671b + 57477a + 1566a + 129168a$ claimed in \cite[Lemma~2]{Mverify}.

\begin{verbatim}
    gap> Sum( mat{ [ 1, 4, 5, 7, 8 ] } );
    [ 196883, 4371, 275, 4371, 275 ]
\end{verbatim}

\section{The permutation character $(1_{2.\B}^{\M})_{2.\B}$}\label{suborbits}

According to \cite[Tables VII, IX]{GMS89},
the restriction of the permutation character $1_{2.\B}^{\M}$ to $2.\B$
decomposes into nine transitive permutation characters $1_U^{2.\B}$,
with the point stabilizers $U$ listed in Table~\ref{suborbitsTable}.

\begin{table}
\caption{Suborbit information}\label{suborbitsTable}
\begin{center}
\begin{tabular}[t]{l|l|r}
    $a c \in$    & $G_{a,c}$  &  $|c^{G_a}|$ \\
  \hline
    \texttt{1A} & $2.\B$            &                    1 \\
    \texttt{2A} & $2^2.{}^2\E_6(2)$ &          27143910000 \\
    \texttt{2B} & $2^{2+22}.\Co_2$  &       11707448673375 \\
    \texttt{3A} & $\Fi_{23}$        &     2031941058560000 \\
    \texttt{3C} & $\Th$             &    91569524834304000 \\
    \texttt{4A} & $2^{1+22}.\McL$   &  1102935324621312000 \\
    \texttt{4B} & $2.\F_4(2)$       &  1254793905192960000 \\
    \texttt{5A} & $\HN$             & 30434513446055706624 \\
    \texttt{6A} & $2.\Fi_{22}$      & 64353605265653760000 \\
  \hline
\end{tabular}
\end{center}
\end{table}

Here $a$ denotes the central involution in $2.\B$,
the action is that on the $\M$-conjugacy class of $a$,
and $c \in a^{\M}$ is a representative of the orbit in question.

In this section, we compute the nine characters $1_U^{2.\B}$,
where $U$ is one of the above point stabilizers $G_{a,c}$.
Note that $a \in G_{a,c}$ holds
(and thus the character is an inflated character of $\B$)
if and only if $a$ and $c$ commute;
this happens exactly for the first three orbits.

All subgroups $U$ except $2^{2+22}.\Co_2$ and $2^{1+22}.\McL$
are {\ATLAS} groups whose character tables have been verified.
The subgroup $2^{2+22}.\Co_2$ is the preimage of a maximal subgroup
$2^{1+22}.\Co_2$ of $\B$ under the natural epimorphism from $2.\B$,
and the computation/verification of the character table of $2^{1+22}.\Co_2$
has been described in \cite{BMverify}.
It will turn out that we do not need the character table of $2^{1+22}.\McL$.

The nine characters will be stored in the variables
\texttt{pi1}, \texttt{pi2}, \ldots, \texttt{pi9}.

For $U = 2.\B$,
we have $1_U^{2.\B} = 1_{2.\B}$.

\begin{verbatim}
    gap> pi1:= TrivialCharacter( table2B );;
\end{verbatim}

For $U = 2^2.{}^2\E_6(2)$,
the character $1_U^{2.\B}$ is the inflation of $1_{\overline{U}}^{\B}$
from $\B$ to $2.\B$,
for $\overline{U} = U / \langle a \rangle = 2.{}^2\E_6(2)$.
(Note that the class fusion from $\overline{U}$ to $\B$ is not uniquely
determined by the character tables of the two groups,
but the permutation character is unique.)

\begin{verbatim}
    gap> tableB:= CharacterTable( "B" );;   
    gap> tableUbar:= CharacterTable( "2.2E6(2)" );;
    gap> fus:= PossibleClassFusions( tableUbar, tableB );;
    gap> pi:= Set( fus,
    >              map -> InducedClassFunctionsByFusionMap( tableUbar, tableB,
    >                         [ TrivialCharacter( tableUbar ) ], map )[1] );;
    gap> Length( pi );
    1
    gap> pi2:= Inflated( tableB, table2B, pi )[1];;
    gap> mult:= List( Irr( table2B ),
    >                 chi -> ScalarProduct( table2B, chi, pi2 ) );;
    gap> Maximum( mult );
    1
    gap> Positions( mult, 1 );
    [ 1, 2, 3, 5, 7, 13, 15, 17 ]
\end{verbatim}

For $U = 2^{2+22}.\Co_2$,
the character $1_U^{2.\B}$ is the inflation of $1_{\overline{U}}^{\B}$
from $\B$ to $2.\B$,
for $\overline{U} = U / \langle a \rangle = 2^{1+22}.\Co_2$,
a maximal subgroup of $\B$.

\begin{verbatim}
    gap> tableUbar:= CharacterTable( "BN2B" );
    CharacterTable( "2^(1+22).Co2" )
    gap> fus:= PossibleClassFusions( tableUbar, tableB );;
    gap> pi:= Set( fus,
    >              map -> InducedClassFunctionsByFusionMap( tableUbar, tableB,
    >                         [ TrivialCharacter( tableUbar ) ], map )[1] );;
    gap> Length( pi );
    1
    gap> pi3:= Inflated( tableB, table2B, pi )[1];;
    gap> mult:= List( Irr( table2B ),
    >                 chi -> ScalarProduct( table2B, chi, pi3 ) );;
    gap> Maximum( mult );
    1
    gap> Positions( mult, 1 );
    [ 1, 3, 5, 8, 13, 15, 28, 30, 37, 40 ]
\end{verbatim}

Next we consider $U = \Fi_{23}$.

\begin{verbatim}
    gap> tableU:= CharacterTable( "Fi23" );;
    gap> fus:= PossibleClassFusions( tableU, table2B );;
    gap> pi:= Set( fus,
    >              map -> InducedClassFunctionsByFusionMap( tableU, table2B,
    >                         [ TrivialCharacter( tableU ) ], map )[1] );;
    gap> Length( pi );
    1
    gap> pi4:= pi[1];;
    gap> mult:= List( Irr( table2B ),
    >                 chi -> ScalarProduct( table2B, chi, pi4 ) );;
    gap> Maximum( mult );
    1
    gap> Positions( mult, 1 );
    [ 1, 2, 3, 5, 7, 8, 9, 12, 13, 15, 17, 23, 27, 30, 32, 40, 41, 54, 
      63, 68, 77, 81, 83, 185, 186, 187, 188, 189, 194, 195, 196, 203, 
      208, 220 ]
\end{verbatim}

Next we consider $U = \Th$.

\begin{verbatim}
    gap> tableU:= CharacterTable( "Th" );;
    gap> fus:= PossibleClassFusions( tableU, table2B );;
    gap> pi:= Set( fus,
    >              map -> InducedClassFunctionsByFusionMap( tableU, table2B,
    >                         [ TrivialCharacter( tableU ) ], map )[1] );;
    gap> Length( pi );
    1
    gap> pi5:= pi[1];;
    gap> mult:= List( Irr( table2B ),
    >                 chi -> ScalarProduct( table2B, chi, pi5 ) );;
    gap> Maximum( mult );
    2
    gap> Positions( mult, 1 );
    [ 1, 3, 7, 8, 12, 13, 16, 19, 27, 28, 34, 38, 41, 57, 68, 70, 77, 78, 
      85, 89, 113, 114, 116, 129, 133, 142, 143, 145, 155, 156, 185, 187, 
      188, 193, 195, 196, 201, 208, 216, 219, 225, 232, 233, 235, 236, 
      237, 242 ]
    gap> Positions( mult, 2 );
    [ 62 ]
\end{verbatim}

For $U = 2^{1+22}.\McL$, we carry out the computations described in
\cite[Section ``A permutation character of $2.\B$'']{ctblpope}.
We know that $U$ is a subgroup of $2^{2+22}.\Co_2$,
and that $\langle U, a \rangle$ has the structure $2^{2+22}.\McL$.

As a first step, we induce the trivial character of $\langle U, a \rangle$
to $2.\B$,
which can be performed by inducing the trivial character of $\McL$ to $\Co_2$,
then to inflate this character to $2^{1+22}.\Co_2$,
then to induce this character to $\B$,
and then to inflate this character to $2.\B$,

\begin{verbatim}
    gap> mcl:= CharacterTable( "McL" );;
    gap> co2:= CharacterTable( "Co2" );;
    gap> fus:= PossibleClassFusions( mcl, co2 );;       
    gap> Length( fus );
    4
    gap> ind:= Set( fus, map -> InducedClassFunctionsByFusionMap( mcl, co2,     
    >                               [ TrivialCharacter( mcl ) ], map )[1] );;
    gap> Length( ind );
    1
    gap> bm2:= CharacterTable( "BM2" );
    CharacterTable( "2^(1+22).Co2" )
    gap> infl:= Inflated( co2, bm2, ind );;
    gap> ind:= Induced( bm2, tableB, infl );;
    gap> infl:= Inflated( tableB, table2B, ind )[1];;
\end{verbatim}

As a second step,
we compute $1_U^{2.\B}$ with the {\GAP} function \texttt{PermChars},
using that we can speed up these computations by prescribing
the permutation character induced from the closure of $U$ with
the normal subgroup $\langle a \rangle$ of $2.\B$.

(We are lucky:
There is a unique solution, and its computation is quite fast.)

\begin{verbatim}
    gap> centre:= ClassPositionsOfCentre( table2B );
    [ 1, 2 ]
    gap> pi:= PermChars( table2B, rec( torso:= [ 2 * infl[1], 0 ],
    >                             normalsubgroup:= centre,
    >                             nonfaithful:= infl ) );;
    gap> Length( pi );
    1
    gap> pi6:= pi[1];;
    gap> List( Irr( table2B ), chi -> ScalarProduct( table2B, chi, pi6 ) );
    [ 1, 1, 2, 1, 2, 0, 2, 3, 2, 0, 0, 1, 4, 1, 2, 0, 3, 2, 0, 2, 0, 0, 
      2, 2, 0, 0, 2, 3, 1, 5, 0, 4, 3, 2, 0, 0, 3, 2, 0, 6, 4, 0, 1, 1, 
      0, 0, 0, 0, 3, 0, 1, 0, 0, 5, 0, 5, 2, 0, 0, 2, 0, 0, 4, 1, 0, 2, 
      0, 4, 2, 4, 4, 3, 0, 2, 4, 2, 4, 0, 3, 0, 3, 2, 5, 0, 1, 0, 3, 1, 
      0, 1, 1, 2, 5, 3, 1, 1, 4, 5, 1, 1, 0, 3, 0, 0, 3, 2, 1, 1, 2, 1, 
      1, 4, 0, 3, 2, 3, 1, 3, 0, 1, 3, 0, 2, 2, 1, 3, 3, 0, 0, 2, 0, 0, 
      0, 0, 3, 0, 3, 3, 3, 1, 0, 3, 0, 4, 0, 1, 0, 0, 2, 0, 0, 2, 0, 0, 
      2, 1, 1, 0, 0, 0, 0, 1, 2, 1, 1, 1, 0, 1, 1, 1, 1, 1, 1, 0, 2, 1, 
      1, 3, 3, 0, 0, 0, 1, 1, 1, 1, 2, 3, 2, 0, 0, 2, 2, 4, 3, 5, 2, 4, 
      0, 0, 0, 0, 5, 2, 0, 0, 0, 1, 1, 0, 0, 0, 0, 0, 0, 7, 0, 0, 1, 7, 
      7, 0, 0, 0, 1, 6, 4, 5, 0, 0, 3, 0, 0, 0, 0, 0, 4, 1, 1, 3, 8, 3, 
      2, 2, 5, 0, 1 ]
\end{verbatim}

Next we consider $U = 2.\F_4(2)$.
We know that $U$ does not contain the central involution of $2.\B$.

\begin{verbatim}
    gap> tableU:= CharacterTable( "2.F4(2)" );;
    gap> fus:= PossibleClassFusions( tableU, table2B );;
    gap> pi:= Set( fus, map -> InducedClassFunctionsByFusionMap( tableU, table2B,
    >             [ TrivialCharacter( tableU ) ], map )[1] );;
    gap> Length( pi );
    2
    gap> pi:= Filtered( pi, x -> ClassPositionsOfKernel( x ) = [ 1 ] );;
    gap> Length( pi );
    1
    gap> pi7:= pi[1];;
    gap> List( Irr( table2B ), chi -> ScalarProduct( table2B, chi, pi7 ) );
    [ 1, 1, 2, 0, 2, 0, 2, 2, 1, 0, 0, 2, 4, 1, 3, 0, 2, 1, 0, 0, 0, 0, 
      2, 1, 0, 0, 2, 2, 1, 4, 0, 2, 1, 2, 0, 0, 3, 2, 0, 4, 4, 0, 0, 0, 
      0, 0, 1, 0, 0, 0, 1, 0, 0, 2, 1, 3, 3, 0, 0, 3, 0, 1, 4, 0, 0, 3, 
      0, 6, 0, 3, 2, 0, 0, 1, 4, 1, 4, 2, 6, 1, 4, 0, 4, 0, 1, 1, 2, 0, 
      0, 3, 2, 1, 3, 2, 0, 0, 4, 5, 3, 1, 0, 3, 0, 0, 1, 1, 2, 0, 0, 2, 
      0, 2, 0, 3, 3, 3, 0, 4, 1, 0, 4, 1, 1, 1, 1, 1, 2, 1, 1, 2, 3, 0, 
      0, 2, 2, 0, 5, 5, 3, 0, 1, 5, 1, 4, 0, 1, 0, 1, 1, 0, 0, 3, 1, 0, 
      2, 3, 1, 0, 2, 0, 0, 2, 1, 0, 1, 0, 0, 1, 0, 0, 0, 1, 0, 0, 1, 2, 
      1, 4, 4, 0, 0, 0, 3, 1, 1, 1, 2, 2, 2, 0, 0, 1, 2, 3, 3, 3, 1, 2, 
      0, 0, 1, 1, 4, 2, 0, 0, 0, 3, 2, 0, 0, 0, 0, 0, 0, 4, 0, 0, 1, 5, 
      5, 0, 1, 1, 2, 2, 4, 4, 0, 0, 3, 1, 1, 1, 0, 0, 4, 1, 1, 5, 7, 3, 
      2, 5, 5, 0, 1 ]
\end{verbatim}

Next we consider $U = \HN$.

\begin{verbatim}
    gap> tableU:= CharacterTable( "HN" );;
    gap> fus:= PossibleClassFusions( tableU, table2B );;
    gap> pi:= Set( fus, map -> InducedClassFunctionsByFusionMap( tableU, table2B,
    >             [ TrivialCharacter( tableU ) ], map )[1] );;
    gap> Length( pi );
    1
    gap> pi8:= pi[1];;
    gap> List( Irr( table2B ), chi -> ScalarProduct( table2B, chi, pi8 ) );
    [ 1, 1, 2, 1, 2, 0, 3, 4, 2, 1, 1, 4, 4, 2, 1, 1, 3, 3, 1, 3, 0, 0, 
      5, 3, 0, 0, 6, 4, 5, 6, 1, 7, 4, 7, 0, 0, 3, 8, 2, 6, 11, 2, 5, 5, 
      0, 0, 2, 1, 3, 4, 7, 0, 0, 7, 3, 9, 5, 0, 0, 6, 4, 2, 13, 6, 0, 4, 
      4, 12, 11, 16, 9, 7, 3, 11, 13, 12, 20, 5, 10, 6, 11, 13, 17, 4, 
      10, 7, 19, 7, 7, 8, 10, 14, 18, 19, 5, 10, 12, 23, 7, 12, 6, 24, 6, 
      4, 17, 16, 8, 9, 17, 11, 12, 23, 8, 24, 18, 26, 21, 29, 10, 18, 31, 
      10, 24, 21, 17, 27, 35, 13, 14, 29, 19, 12, 7, 18, 26, 15, 34, 34, 
      35, 20, 14, 36, 14, 39, 8, 29, 24, 15, 40, 13, 9, 38, 24, 17, 35, 
      32, 26, 26, 24, 22, 17, 31, 39, 29, 30, 30, 19, 44, 37, 37, 28, 30, 
      31, 29, 42, 40, 40, 56, 56, 30, 30, 42, 50, 47, 2, 2, 4, 6, 4, 0, 
      0, 4, 6, 10, 10, 12, 8, 12, 0, 0, 2, 4, 16, 10, 0, 0, 2, 12, 10, 0, 
      0, 0, 0, 0, 0, 28, 0, 0, 14, 34, 40, 2, 10, 10, 22, 40, 44, 44, 8, 
      8, 36, 14, 14, 16, 8, 8, 46, 28, 28, 58, 90, 72, 70, 92, 104, 56, 
      90 ]
\end{verbatim}

Finally, we consider $U = 2.\Fi_{22}$.
There are two candidates for the permutation character $(1_U)^{2.\B}$,
according to the possible class fusions.
One of the two characters is zero on the class of the central involution
of $2.\B$, the other is not.
We know that $U$ does not contain the central involution of $2.\B$,
hence we can decide which character is correct.

\begin{verbatim}
    gap> tableU:= CharacterTable( "2.Fi22" );;
    gap> fus:= PossibleClassFusions( tableU, table2B );;
    gap> pi:= Set( fus, map -> InducedClassFunctionsByFusionMap( tableU, table2B,
    >             [ TrivialCharacter( tableU ) ], map )[1] );;
    gap> Length( pi );
    2
    gap> pi:= Filtered( pi, x -> ClassPositionsOfKernel( x ) = [ 1 ] );;
    gap> Length( pi );
    1
    gap> pi9:= pi[1];;
    gap> List( Irr( table2B ), chi -> ScalarProduct( table2B, chi, pi9 ) );
    [ 1, 2, 3, 1, 4, 1, 5, 5, 5, 1, 1, 5, 8, 4, 4, 1, 7, 6, 0, 5, 0, 0, 
      10, 7, 0, 0, 10, 6, 6, 13, 3, 14, 10, 11, 0, 0, 5, 11, 2, 14, 19, 
      6, 6, 5, 0, 0, 0, 3, 6, 7, 11, 0, 0, 17, 2, 20, 9, 0, 0, 12, 8, 1, 
      23, 11, 1, 8, 7, 23, 18, 27, 18, 12, 7, 22, 29, 21, 34, 6, 22, 7, 
      22, 18, 33, 3, 19, 10, 34, 12, 12, 15, 17, 28, 34, 34, 7, 20, 26, 
      40, 15, 25, 3, 40, 9, 6, 34, 25, 18, 21, 30, 21, 18, 43, 12, 45, 
      39, 49, 38, 51, 18, 32, 63, 19, 42, 41, 33, 48, 64, 27, 29, 52, 38, 
      29, 19, 40, 47, 31, 69, 69, 65, 42, 35, 68, 27, 73, 20, 53, 46, 38, 
      75, 29, 24, 72, 50, 41, 72, 68, 58, 52, 54, 50, 44, 64, 75, 58, 69, 
      65, 49, 85, 75, 75, 63, 68, 65, 63, 90, 87, 83, 118, 118, 74, 71, 
      90, 109, 109, 2, 3, 6, 9, 8, 0, 0, 7, 10, 18, 16, 22, 12, 23, 0, 0, 
      2, 6, 28, 19, 0, 0, 5, 16, 18, 0, 0, 0, 0, 0, 0, 52, 1, 1, 26, 59, 
      76, 11, 18, 18, 39, 77, 80, 77, 22, 22, 66, 27, 27, 33, 20, 20, 87, 
      60, 60, 103, 175, 148, 152, 187, 215, 140, 201 ]
\end{verbatim}

Now we can form the restriction of $(1_{2.\B})^{\M}$ to $2.\B$.

\begin{verbatim}
    gap> constit:= [ pi1, pi2, pi3, pi4, pi5, pi6, pi7, pi8, pi9 ];;
    gap> pi:= Sum( constit );;
\end{verbatim}

\section{The conjugacy classes of $\M$}\label{Mclasses}

\subsection{Our strategy to describe the conjugacy classes of $\M$}%
\label{strategy_classes}

We know the order of $\M$ and its prime divisors.
Let us check whether this fits to our data computed up to now.

\begin{verbatim}
    gap> sizeM:= pi[1] * Size( table2B );
    808017424794512875886459904961710757005754368000000000
    gap> StringPP( sizeM );
    "2^46*3^20*5^9*7^6*11^2*13^3*17*19*23*29*31*41*47*59*71"
    gap> sizeM = Size( CharacterTable( "M" ) );
    true
\end{verbatim}

For each prime $p$ dividing $|\M|$,
we classify the conjugacy classes of elements of order $p$ in $\M$
and use the facts that for each such class representative $x$,
the classes of roots of $x$ in the centralizer/normalizer of $x$
are in bijection with the corresponding classes in $\M$,
and that this bijection respects centralizer orders.

%

For each element $x \in \M$ of order $p \in \{ 2, 3, 5 \}$,
we will use the character table of $N_{\M}(\langle x \rangle)$
to establish $\M$-conjugacy classes of roots of $x$.
In order not to count the same class several times,
we proceed by increasing $p$,
and collect only those classes of roots of $x$ for which $p$ is the smallest
prime divisor of the element order.

For elements $x \in \M$ of prime order $p > 5$,
it is not necessary to use the character table of $N_{\M}(\langle x \rangle)$;
we will use the permutation character values $(1_{2.\B})^{\M}(x)$
and ad hoc arguments.

\subsection{Utility functions}\label{functions}

During the process of finding the conjugacy classes of $\M$,
we record our knowledge about the character table of $\M$
in a global {\GAP} variable \texttt{head},
which is a record with the following components.

\begin{description}
\item[\texttt{Size}] \hfill \\
    the group order $|\M|$,

\item[\texttt{SizesCentralizers}] \hfill \\
    the list of centralizer orders of the conjugacy classes
    established up to now,

\item[\texttt{OrdersClassRepresentatives}] \hfill \\
    the list of corresponding representative orders,

\item[\texttt{fusions}] \hfill \\
    a list that collects the currently known partial class fusions into $\M$;
    each entry is a record with the components
    \texttt{subtable} (the character table of the subgroup)
    and \texttt{map} (the list of known images;
    unknown positions are unbound).
\end{description}

We initialize this variable, using the group order $\M$
and that there is an identity element.

\begin{verbatim}
    gap> head:= rec( Size:= sizeM,
    >                SizesCentralizers:= [ sizeM ],
    >                OrdersClassRepresentatives:= [ 1 ],
    >                fusions:= [],
    >              );;
\end{verbatim}

The function \texttt{ExtendTableHeadByRootClasses} takes
the object \texttt{head},
the character table \texttt{s} of a subgroup $H$ of $\M$,
and an integer \texttt{pos} as its arguments,
where it is assumed that the \texttt{pos}-th class of \texttt{s}
contains an element $x$ of prime order $p$
such that $N_{\M}(\langle x \rangle) = H$ holds
and such that \texttt{head} contains information only about
those classes of $\M$ whose elements have order divisible by a prime
that is smaller than $p$.

\begin{verbatim}
    gap> ExtendTableHeadByRootClasses:= function( head, s, pos )
    >    local fus, orders, p, cents, oldnumber, i, ord;
    > 
    >    # Initialize the fusion information.
    >    fus:= rec( subtable:= s, map:= [ 1 ] );
    >    Add( head.fusions, fus );
    > 
    >    # Compute the positions of root classes of 'pos'.
    >    orders:= OrdersClassRepresentatives( s );
    >    p:= orders[ pos ];
    >    cents:= SizesCentralizers( s );
    >    oldnumber:= Length( head.OrdersClassRepresentatives );
    > 
    >    # Run over the classes of 's'
    >    # are already contained in head
    >    for i in [ 1 .. NrConjugacyClasses( s ) ] do
    >      ord:= orders[i];
    >      if ord mod p = 0 and
    >         Minimum( PrimeDivisors( ord ) ) = p and
    >         PowerMap( s, ord / p, i ) = pos then
    >        # Class 'i' is a root class of 'pos' and is new in 'head'.
    >        Add( head.SizesCentralizers, cents[i] );
    >        Add( head.OrdersClassRepresentatives, orders[i] );
    >        fus.map[i]:= Length( head.SizesCentralizers );
    >      fi;
    >    od;
    > 
    >    Print( "#I  after ", Identifier( s ), ": found ",
    >           Length( head.OrdersClassRepresentatives ) - oldnumber,
    >           " classes, now have ",
    >           Length( head.OrdersClassRepresentatives ), "\n" );
    >    end;;
\end{verbatim}

In several cases, we will establish a conjugacy class $g^{\M}$ without
knowing the character table of a suitable subgroup of $\M$ to which
\texttt{ExtendTableHeadByRootClasses} can be applied, where $g$ is among
the root classes.
That is, we may know just element order \texttt{s}
and centralizer order \texttt{cent}.

We are a bit better off if we know the character table \texttt{s}
of a subgroup of $\M$ and the list \texttt{poss} of all those classes
in this table which fuse to the class $g^{\M}$, because then we can
store this information in the partial class fusion from \texttt{s}
that is stored in \texttt{head}.

\begin{verbatim}
    gap> ExtendTableHeadByCentralizerOrder:= function( head, s, cent, poss )
    >    local ord, fus, i;
    > 
    >    if IsCharacterTable( s ) then
    >      ord:= Set( OrdersClassRepresentatives( s ){ poss } );
    >      if Length( ord ) <> 1 then
    >        Error( "classes cannot fuse" );
    >      fi;
    >      ord:= ord[1];
    >    elif IsInt( s ) then
    >      ord:= s;
    >    fi;
    >    Add( head.SizesCentralizers, cent );
    >    Add( head.OrdersClassRepresentatives, ord );
    > 
    >    Print( "#I  after order ", ord, " element" );
    >    if IsCharacterTable( s ) then
    >      # extend the stored fusion from s
    >      fus:= First( head.fusions,
    >                   r -> Identifier( r.subtable ) = Identifier( s ) );
    >      for i in poss do
    >        fus.map[i]:= Length( head.SizesCentralizers );
    >      od;
    >      Print( " from ", Identifier( s ) );
    >    fi;
    >    Print( ": have ",
    >           Length( head.OrdersClassRepresentatives ), " classes\n" );
    >    end;;
\end{verbatim}

The permutation character $1_H^G$, where $H \leq G$ are two groups,
has the property $1_H^G(g) = |C_G(g)| \cdot |g^G \cap H| / |H|$.
For $g \in H$,
this implies that $|C_G(g)| = 1_H^G(g) \cdot |H| / |g^G \cap H|$
can be computed from the character $(1_H^G)_H$ and the class lengths in $H$,
provided that we know which classes of $H$ fuse into $g^G$.
The function \texttt{ExtendTableHeadByPermCharValue} extends the information
in \texttt{head} by the data for the class $g^{\M}$,
where \texttt{s} is the character table of $H$,
\texttt{pi\_rest\_to\_s} is $(1_H^G)_H$,
and \texttt{poss} is the list of positions of those classes in \texttt{s}
that fuse to $g^{\M}$.

\begin{verbatim}
    gap> ExtendTableHeadByPermCharValue:= function( head, s, pi_rest_to_s, poss )
    >    local pival, cent;
    > 
    >    pival:= Set( pi_rest_to_s{ poss } );
    >    if Length( pival ) <> 1 then
    >      Error( "classes cannot fuse" );
    >    fi;
    > 
    >    cent:= pival[1] * Size( s ) / Sum( SizesConjugacyClasses( s ){ poss } );
    >    ExtendTableHeadByCentralizerOrder( head, s, cent, poss );
    >    end;;
\end{verbatim}

\subsection{Classes of elements of even order}\label{elements_2}

By~\cite{Mverify},
we know that $\M$ has exactly two conjugacy classes of involutions,
and that the involution centralizers have the structures
$2.\B$ (for the class \texttt{2A}) and
$2^{1+24}_+.\Co_1$ (for the class \texttt{2B}), respectively.

Moreover, the character tables of these subgroups that are
stored in the {\GAP} Character Table Library are correct.
For $2.\B$, this follows from the correctness of the character table of $\B$
as shown in \cite{BMverify} and the computations
in \cite{AtlasVerifyLargeArxiv}.
For $2^{1+24}_+.\Co_1$, the recomputation of the character table is described
in Section~\ref{sect:table_c2b}.

Thus we can determine the $\M$-conjugacy classes of elements of even order
as follows.

\begin{verbatim}
    gap> s:= CharacterTable( "2.B" );;
    gap> ClassPositionsOfCentre( s );
    [ 1, 2 ]
    gap> ExtendTableHeadByRootClasses( head, s, 2 );
    #I  after 2.B: found 42 classes, now have 43
    gap> s:= CharacterTable( "MN2B" );;
    gap> ClassPositionsOfCentre( s );
    [ 1, 2 ]
    gap> ExtendTableHeadByRootClasses( head, s, 2 );
    #I  after 2^1+24.Co1: found 91 classes, now have 134
\end{verbatim}

\subsection{Classes of elements of order divisible by $3$}\label{elements_3}

We know that $\M$ has exactly three conjugacy classes of elements
of order $3$,
and that their normalizers have the structures
$3.\Fi_{24}$ (for the class \texttt{3A}),
$3^{1+12}_+.2.\Suz.2$ (for the class \texttt{3B}),
and $S_3 \times \Th$ (for the class \texttt{3C}), respectively.

Moreover,
the {\GAP} character tables of $3.\Fi_{24}$ and $\Th$ are {\ATLAS} tables
and have been verified, see~\cite{BMO17}.

We determine the $\M$-conjugacy classes of elements of odd order
that are roots of \texttt{3A} or \texttt{3C} elements, as follows.

\begin{verbatim}
    gap> s:= CharacterTable( "3.Fi24" );;
    gap> ClassPositionsOfPCore( s, 3 );
    [ 1, 2 ]
    gap> ExtendTableHeadByRootClasses( head, s, 2 );
    #I  after 3.F3+.2: found 12 classes, now have 146
    gap> s:= CharacterTableDirectProduct( CharacterTable( "Th" ),
    >                                     CharacterTable( "Symmetric", 3 ) );;
    gap> ClassPositionsOfPCore( s, 3 );
    [ 1, 3 ]
    gap> ExtendTableHeadByRootClasses( head, s, 3 );
    #I  after ThxSym(3): found 7 classes, now have 153
\end{verbatim}

The situation with the \texttt{3B} normalizer is more involved.
Section~\ref{sect:norm3B} documents the construction of the character table
of a downward extension of the structure $3^{1+12}_+:6.\Suz.2$
of the \texttt{3B} normalizer, and gives two
candidates for the character table of the \texttt{3B} normalizer.

It will turn out that each of these candidates leads to ``the same''
root classes,
in the sense that the number of these classes, their element orders,
and their centralizer orders are equal.
Note that the $3$-core of $H = 3^{1+12}_+:6.\Suz.2$ has the structure
$X \times N$,
where $X$ has order $3$
and $N \cong 3^{1+12}_+$ such that $H / N \cong 6.\Suz.2$ holds.
We are interested in the two ``diagonal'' factors, that is,
the factors of $H$ by the one of the two normal subgroups of order $3$ in $H$
that are not equal to $X$ or $Z(N)$.
(See the picture in Section~\ref{sect:norm3B} for the details.)

First we exclude the normal subgroup of order $3$ that is contained in the
unique normal subgroup $N$ of order $3^{13}$.

\begin{verbatim}
    gap> exts:= CharacterTable( "3^(1+12):6.Suz.2" );;
    gap> kernels:= Positions( SizesConjugacyClasses( exts ), 2 );
    [ 2, 18, 19, 20 ]
    gap> order3_13:= Filtered( ClassPositionsOfNormalSubgroups( exts ),
    >        l -> Sum( SizesConjugacyClasses( exts ){ l } ) = 3^13 );
    [ [ 1 .. 4 ] ]
    gap> kernels:= Difference( kernels, order3_13[1] );
    [ 18, 19, 20 ]
\end{verbatim}

The classes in the subgroup $X$ can be identified by the fact that
exactly one factor of $H$ by a normal subgroup of order $3$ admits a
class fusion from $2.\Suz.2$, and hence this must be the split extension
of $3^{1+12}_+$ with $2.\Suz.2$.


\begin{verbatim}
    gap> facts:= List( kernels, i -> exts / [ 1, i ] );
    [ CharacterTable( "3^(1+12):6.Suz.2/[ 1, 18 ]" ), 
      CharacterTable( "3^(1+12):6.Suz.2/[ 1, 19 ]" ), 
      CharacterTable( "3^(1+12):6.Suz.2/[ 1, 20 ]" ) ]
    gap> f:= CharacterTable( "2.Suz.2" );;
    gap> facts:= Filtered( facts,
    >        x -> Length( PossibleClassFusions( f, x ) ) = 0 );
    [ CharacterTable( "3^(1+12):6.Suz.2/[ 1, 19 ]" ), 
      CharacterTable( "3^(1+12):6.Suz.2/[ 1, 20 ]" ) ]
\end{verbatim}

We compute the root classes for both candidates.
For that,
we first create a copy \texttt{head2} of the information in \texttt{head}.

\begin{verbatim}
    gap> kernels:= List( facts,
    >        f -> Positions( SizesConjugacyClasses( f ), 2 ) );
    [ [ 2 ], [ 2 ] ]
    gap> head2:= StructuralCopy( head );;
    gap> ExtendTableHeadByRootClasses( head, facts[1], 2 );
    #I  after 3^(1+12):6.Suz.2/[ 1, 19 ]: found 12 classes, now have 165
    gap> ExtendTableHeadByRootClasses( head2, facts[2], 2 );
    #I  after 3^(1+12):6.Suz.2/[ 1, 20 ]: found 12 classes, now have 165
\end{verbatim}

We observe that \texttt{head} and \texttt{head2} differ only by the
two character tables in the last fusion record.

\begin{verbatim}
    gap> nams:= RecNames( head );
    [ "Size", "OrdersClassRepresentatives", "SizesCentralizers", 
      "fusions" ]
    gap> ForAll( Difference( nams, [ "fusions" ] ),
    >            nam -> head.( nam ) = head2.( nam ) );
    true
    gap> Length( head.fusions );
    5
    gap> ForAll( [ 1 .. 4 ], i -> head.fusions[i] = head2.fusions[i] );
    true
    gap> head.fusions[5].map = head2.fusions[5].map;
    true
\end{verbatim}

We continue with establishing the conjugacy classes of $\M$.
The question which of the two above candidate tables belongs to a subgroup
of $\M$ will be answered in Section~\ref{sect:natcharM}.

\subsection{Classes of elements of order divisible by $5$}\label{elements_5}

The group $2.\B$ contains two rational conjugacy classes of elements
of order $5$,
with different values in the permutation character $(1_{2.\B})^{\M}$.

\begin{verbatim}
    gap> s:= CharacterTable( "2.B" );;
    gap> pos:= Positions( OrdersClassRepresentatives( s ), 5 );
    [ 23, 25 ]
    gap> pi{ pos };
    [ 1539000, 7875 ]
\end{verbatim}

This establishes two classes \texttt{5A}, \texttt{5B} of conjugacy classes
of elements of order $5$ in $\M$,
with centralizer orders $5 |\HN|$ and $5^7 |2.\J_2|$, respectively.

\begin{verbatim}
    gap> cents:= List( pos,
    >      i -> pi[i] * Size( s ) / SizesConjugacyClasses( s )[i] );
    [ 1365154560000000, 94500000000 ]
    gap> cents = [ 5 * Size( CharacterTable( "HN" ) ),
    >              5^7 * Size( CharacterTable( "2.J2" ) ) ];
    true
\end{verbatim}


By~\cite{Mverify},
we know that $\M$ contains exactly two conjugacy classes of elements
of order $5$,
\texttt{5A} with centralizer $5 \times \HN$ and normalizer
$(D_{10} \times \HN).2$,
and \texttt{5B} with centralizer $5^{1+6}_+.2.\J_2$ and normalizer
$5^{1+6}_+.4.\J_2.2$.

The two classes are rational
because this is the case already for their intersections with $2.\B$.

\begin{verbatim}
    gap> s:= CharacterTable( "MN5A" );
    CharacterTable( "(D10xHN).2" )
    gap> ClassPositionsOfPCore( s, 5 );
    [ 1, 45 ]
    gap> ExtendTableHeadByRootClasses( head, s, 45 );
    #I  after (D10xHN).2: found 5 classes, now have 170
\end{verbatim}

The character table of $5^{1+6}_+.4.\J_2.2$ has been recomputed
with {\MAGMA}, see Section~\ref{sect:table_N5B},
thus we are allowed to use the character table from
the {\GAP} character table library.

\begin{verbatim}
    gap> s:= CharacterTable( "MN5B" );
    CharacterTable( "5^(1+6):2.J2.4" )
    gap> 5core:= ClassPositionsOfPCore( s, 5 );
    [ 1 .. 4 ]
    gap> SizesConjugacyClasses( s ){ 5core };
    [ 1, 4, 37800, 40320 ]
    gap> ExtendTableHeadByRootClasses( head, s, 2 );
    #I  after 5^(1+6):2.J2.4: found 3 classes, now have 173
\end{verbatim}

\subsection{Classes of elements of order divisible by $11$}\label{elements_11}

The group $2.\B$ contains a rational class of elements of order $11$.
The permutation character $(1_{2.\B}^{\M})_{2.\B}$ yields a class of elements
of order $11$ with centralizer order $11 |M_{12}|$ in $\M$.

\begin{verbatim}
    gap> s:= CharacterTable( "2.B" );;
    gap> pos:= Positions( OrdersClassRepresentatives( s ), 11 );
    [ 71 ]
\end{verbatim}

By the arguments in~\cite{Mverify},
$\M$ has no other classes of element order $11$.

\begin{verbatim}
    gap> ExtendTableHeadByPermCharValue( head, s, pi, pos );
    #I  after order 11 element from 2.B: have 174 classes
\end{verbatim}

\subsection{Classes of elements of the orders $17, 19, 23, 31, 47$}%
\label{elements_17}

The elements of the orders $17, 19, 23, 31, 47$ in $\M$ lie in cyclic
Sylow subgroups that appear already in $2.\B$.

The elements of order $17$ and $19$ are rational in $2.\B$
and hence also in $\M$.

\begin{verbatim}
    gap> s:= CharacterTable( "2.B" );;
    gap> pos:= Positions( OrdersClassRepresentatives( s ), 17 );
    [ 118 ]
    gap> ExtendTableHeadByPermCharValue( head, s, pi, pos );
    #I  after order 17 element from 2.B: have 175 classes
    gap> pos:= Positions( OrdersClassRepresentatives( s ), 19 );
    [ 128 ]
    gap> ExtendTableHeadByPermCharValue( head, s, pi, pos );
    #I  after order 19 element from 2.B: have 176 classes
\end{verbatim}

For elements $g$ of order $p \in \{ 23, 31, 47 \}$,
the group $2.\B$ contains exactly two Galois conjugate classes that contain
the nonidentity powers of $g$,
which means that $[N_{2.\B}(\langle g \rangle):C_{2.\B}(g)] = (p-1)/2$ holds.
The equation $|C_{\M}(g)| = |2.\B| \cdot \pi(g) / |g^{\M} \cap 2.\B|$
implies
\[
  |N_{\M}(\langle g \rangle)|
    = [N_{\M}(\langle g \rangle):C_{\M}(g)]
      \cdot |2.\B| \cdot \pi(g) / |g^{\M}|
    = (p-1)/2 \cdot |2.\B| \cdot \pi(g) / |g^{2.\B}|.
\]
Note that either the two classes of elements of order $p$ in $2.\B$
fuse in $\M$ or not;
in the former case,
we have $[N_{\M}(\langle g \rangle):C_{\M}(g)] = p-1$ and
$|g^{\M} \cap 2.\B| = 2 |g^{2.\B}|$,
whereas we have
$[N_{\M}(\langle g \rangle):C_{\M}(g)] = (p-1)/2$ and
$|g^{\M} \cap 2.\B| = |g^{2.\B}|$ in the latter case.
Thus we can compute $|N_{\M}(\langle g \rangle)|$ in each case,
and we can then find arguments why the two Galois conjugate classes
do not fuse.

First we deal with $p = 23$.

\begin{verbatim}
    gap> s:= CharacterTable( "2.B" );
    CharacterTable( "2.B" )
    gap> ord:= OrdersClassRepresentatives( s );;
    gap> classes:= SizesConjugacyClasses( s );;
    gap> p:= 23;;
    gap> pos:= Positions( ord, p );
    [ 147, 149 ]
    gap> n:= (p-1)/2 * Size( s ) * pi[ pos[1] ] / classes[ pos[1] ];
    6072
    gap> Collected( Factors( n ) );
    [ [ 2, 3 ], [ 3, 1 ], [ 11, 1 ], [ 23, 1 ] ]
\end{verbatim}

In order to prove that the two classes of elements of order $23$ in $2.\B$
do not fuse in $\M$, it suffices to show that the centralizer order is
divisible by $2^3$.
We see that this is the case already in the \texttt{2B} centralizer in $\M$.

\begin{verbatim}
    gap> u:= CharacterTable( "MN2B" );
    CharacterTable( "2^1+24.Co1" )
    gap> upos:= Positions( OrdersClassRepresentatives( u ), p );
    [ 289, 294 ]
    gap> SizesCentralizers( u ){ upos } / 2^3;
    [ 23, 23 ]
\end{verbatim}

Thus we have established two classes of element order $23$ in $\M$.

\begin{verbatim}
    gap> ExtendTableHeadByPermCharValue( head, s, pi, pos{ [1] } );
    #I  after order 23 element from 2.B: have 177 classes
    gap> ExtendTableHeadByPermCharValue( head, s, pi, pos{ [2] } );
    #I  after order 23 element from 2.B: have 178 classes
\end{verbatim}

The case $p = 31$ is done analogously.
Here the necessary $2$-part of the centralizer occurs already in $2.\B$.

\begin{verbatim}
    gap> p:= 31;;
    gap> pos:= Positions( OrdersClassRepresentatives( s ), p );
    [ 190, 192 ]
    gap> n:= (p-1)/2 * Size( s ) * pi[ pos[1] ] / classes[ pos[1] ];
    2790
    gap> Collected( Factors( n ) );
    [ [ 2, 1 ], [ 3, 2 ], [ 5, 1 ], [ 31, 1 ] ]
    gap> SizesCentralizers( s ){ pos };
    [ 62, 62 ]
    gap> ExtendTableHeadByPermCharValue( head, s, pi, pos{ [1] } );
    #I  after order 31 element from 2.B: have 179 classes
    gap> ExtendTableHeadByPermCharValue( head, s, pi, pos{ [2] } );
    #I  after order 31 element from 2.B: have 180 classes
\end{verbatim}

Finally, we deal with $p = 47$.

\begin{verbatim}
    gap> p:= 47;;
    gap> pos:= Positions( OrdersClassRepresentatives( s ), p );
    [ 228, 230 ]
    gap> n:= (p-1)/2 * Size( s ) * pi[ pos[1] ] / classes[ pos[1] ];
    2162
    gap> Collected( Factors( n ) );
    [ [ 2, 1 ], [ 23, 1 ], [ 47, 1 ] ]
    gap> SizesCentralizers( s ){ pos };
    [ 94, 94 ]
    gap> ExtendTableHeadByPermCharValue( head, s, pi, pos{ [1] } );
    #I  after order 47 element from 2.B: have 181 classes
    gap> ExtendTableHeadByPermCharValue( head, s, pi, pos{ [2] } );
    #I  after order 47 element from 2.B: have 182 classes
\end{verbatim}

\subsection{Classes of elements of order $13$}\label{elements_13}

The class \texttt{13A} of $\M$ arises from the rational class of elements
of order $13$ in $2.\B$.
We use the permutation character to enter the information about
the class \texttt{13A}.

\begin{verbatim}
    gap> p:= 13;;
    gap> pos:= Positions( OrdersClassRepresentatives( s ), p );
    [ 97 ]
    gap> c:= Size( s ) * pi[ pos[1] ] / classes[ pos[1] ];
    73008
    gap> Factors( c );
    [ 2, 2, 2, 2, 3, 3, 3, 13, 13 ]
    gap> ExtendTableHeadByPermCharValue( head, s, pi, pos );
    #I  after order 13 element from 2.B: have 183 classes
\end{verbatim}

The class \texttt{13B} intersects the \texttt{2B} centralizer.
Here we just know the centralizer order $13^3 \cdot 2^3 \cdot 3$.

\begin{verbatim}
    gap> c2b:= CharacterTable( "MN2B" );;
    gap> pos:= Positions( OrdersClassRepresentatives( c2b ), 13 );
    [ 220 ]
    gap> ExtendTableHeadByCentralizerOrder( head, c2b, 13^3 * 24, pos );
    #I  after order 13 element from 2^1+24.Co1: have 184 classes
\end{verbatim}

\subsection{Classes of elements of order divisible by $29$}\label{elements_29}

The group $3.\Fi_{24}$ contains a rational class of elements of order $29$,
with centralizer order $3 \cdot 29$.

\begin{verbatim}
    gap> u:= CharacterTable( "3.Fi24" );;
    gap> pos:= Positions( OrdersClassRepresentatives( u ), 29 );
    [ 142 ]
    gap> SizesCentralizers( u ){ pos };
    [ 87 ]
\end{verbatim}

The list of classes of $\M$ collected up to now covers all roots of
elements of the orders $2, 3, 5, 11, 13, 17, 19, 23, 31, 47$,
and $29$ occurs as a factor of the centralizer order only for the
classes \texttt{1A}, \texttt{3A}, \texttt{87A}, and \texttt{87B}.

\begin{verbatim}
    gap> poss:= PositionsProperty( head.SizesCentralizers,
    >                              x -> x mod 29 = 0 );
    [ 1, 135, 144, 145 ]
    gap> head.OrdersClassRepresentatives{ poss };
    [ 1, 3, 87, 87 ]
    gap> head.SizesCentralizers{ poss };
    [ 808017424794512875886459904961710757005754368000000000, 
      3765617127571985163878400, 87, 87 ]
\end{verbatim}

Thus the only possible additional prime divisors of the centralizer order
in $\M$ of an element $x$ of order $29$ are $7, 41, 59$, and $71$.

\begin{verbatim}
    gap> candprimes:= Difference( PrimeDivisors( head.Size ),
    >                     [ 2, 3, 5, 11, 13, 17, 19, 23, 29, 31, 47 ] );
    [ 7, 41, 59, 71 ]
\end{verbatim}

The centralizer order of $x$ has the form
$3 \cdot 29 \cdot 7^i \cdot 41^j \cdot 59^k \cdot 71^l$,
with $ 0 \leq i \leq 6$ and $j, k, l \in \{ 0, 1 \}$.

\begin{verbatim}
    gap> parts:= Filtered( Collected( Factors( head.Size ) ),
    >                      x -> x[1] in candprimes );
    [ [ 7, 6 ], [ 41, 1 ], [ 59, 1 ], [ 71, 1 ] ]
    gap> poss:= List( parts, l -> List( [ 0 .. l[2] ], i -> l[1]^i ) );;
    gap> cart:= Cartesian( poss );;
    gap> possord:= 3 * 29 * List( cart, Product );;
\end{verbatim}

Only $3 \cdot 29$ and $3 \cdot 29 \cdot 59$ satisfy Sylow's theorem,
that is, $|\M| / |N_{\M}(\langle x \rangle)| \equiv 1 \pmod{29}$.
Note that we have $[N_{\M}(\langle x \rangle):C_{\M}(x)] = 28$.

\begin{verbatim}
    gap> good:= Filtered( possord,
    >                     x -> ( head.Size / ( 28 * x ) ) mod 29 = 1 );
    [ 87, 5133 ]
    gap> List( good, Factors );
    [ [ 3, 29 ], [ 3, 29, 59 ] ]
\end{verbatim}

Now we can exclude the possible centralizer order $3 \cdot 29 \cdot 59$
by the fact that the Sylow $59$ subgroup would be normal
and thus would be normalized and hence centralized by an element of
order $3$, a contradiction.

\begin{verbatim}
    gap> Filtered( DivisorsInt( 5133 ), x -> x mod 59 = 1 );
    [ 1 ]
\end{verbatim}

Thus we have established a rational class of elements of order $29$,
with centralizer of order $3 \cdot 29$.

\begin{verbatim}
    gap> ExtendTableHeadByCentralizerOrder( head, u, 3 * 29, pos );
    #I  after order 29 element from 3.F3+.2: have 185 classes
\end{verbatim}

\subsection{Classes of elements of order divisible by $41$}\label{elements_41}

We assume that $\M$ contains a subgroup of the structure $3^8.O_8^-(3)$.
The fact that an element $x$ of order $41$ in $\M$ is normalized
by an element of order $4$ can be read off from the factor group
$O_8^-(3)$.

\begin{verbatim}
    gap> t:= CharacterTable( "O8-(3)" );
    CharacterTable( "O8-(3)" )
    gap> Length( Positions( OrdersClassRepresentatives( t ), 41 ) );
    10
\end{verbatim}

By the above arguments, the only possible odd prime divisors of
$|N_{\M}(\langle x \rangle)|/41$ are $5, 7, 59, 71$,
where $5$ cannot divide the centralizer order.
As in the case of $p = 29$, we apply Sylow's theorem,
and get
$|N_{\M}(\langle x \rangle)| \in
\{ 2^3 \cdot 5 \cdot 41, 2^3 \cdot 7 \cdot 41 \cdot 71 \}$.

\begin{verbatim}
    gap> possord:= 2^2 * 41 * DivisorsInt( 2 * 5 * 7^6 * 59 * 71 );;
    gap> good:= Filtered( possord,
    >                     x -> ( head.Size / x ) mod 41 = 1 );
    [ 1640, 163016 ]
    gap> List( good, Factors );
    [ [ 2, 2, 2, 5, 41 ], [ 2, 2, 2, 7, 41, 71 ] ]
\end{verbatim}

Suppose that $71$ divides $|N_{\M}(\langle x \rangle)|$.
Then the $71$ Sylow subgroup of $N_{\M}(\langle x \rangle)$ is normal
thus normalized by an element of order $8$,
and thus centralized by an involution, a contradiction.

\begin{verbatim}
    gap> Filtered( DivisorsInt( good[2] ), x -> x mod 71 = 1 );
    [ 1 ]
\end{verbatim}

Thus we have established a rational class of self-centralizing elements
of order $41$.

\begin{verbatim}
    gap> ExtendTableHeadByCentralizerOrder( head, 41, 41, fail );
    #I  after order 41 element: have 186 classes
\end{verbatim}

\subsection{Classes of elements of order divisible by $59$}\label{elements_59}

By the above arguments,
the normalizer order of an element of order $59$ divides
$58 \cdot 7^6 \cdot 59 \cdot 71$.
Sylow's theorem admits just the normalizer order $59 \cdot 29$.

\begin{verbatim}
    gap> possord:= 59 * DivisorsInt( 58*7^6*71 );;
    gap> good:= Filtered( possord,
    >                     x -> ( head.Size / x ) mod 59 = 1 );
    [ 1711 ]
    gap> List( good, Factors );
    [ [ 29, 59 ] ]
\end{verbatim}

Thus we have established a pair of Galois conjugate classes of
self-centralizing elements of order $59$.

\begin{verbatim}
    gap> ExtendTableHeadByCentralizerOrder( head, 59, 59, fail );
    #I  after order 59 element: have 187 classes
    gap> ExtendTableHeadByCentralizerOrder( head, 59, 59, fail );
    #I  after order 59 element: have 188 classes
\end{verbatim}

\subsection{Classes of elements of order divisible by $71$}\label{elements_71}

By the above arguments, 
the normalizer order of an element of order $71$ divides
$70 \cdot 7^5 \cdot 71$.
Sylow's theorem admits just the normalizer order $71 \cdot 35$.

\begin{verbatim}
    gap> possord:= 71 * DivisorsInt( 70*7^5 );;
    gap> good:= Filtered( possord,
    >                     x -> ( head.Size / x ) mod 71 = 1 );
    [ 2485 ]
    gap> List( good, Factors );
    [ [ 5, 7, 71 ] ]
\end{verbatim}

Thus we have established a pair of Galois conjugate classes of
self-centralizing elements of order $71$.

\begin{verbatim}
    gap> ExtendTableHeadByCentralizerOrder( head, 71, 71, fail );
    #I  after order 71 element: have 189 classes
    gap> ExtendTableHeadByCentralizerOrder( head, 71, 71, fail );
    #I  after order 71 element: have 190 classes
\end{verbatim}

\subsection{Classes of elements of order divisible by $7$}\label{elements_7}

The subgroup $2.\B$ yields a rational class \texttt{7A}
with centralizer order $7 \cdot |\He|$.

\begin{verbatim}
    gap> s:= CharacterTable( "2.B" );;
    gap> pos:= Positions( OrdersClassRepresentatives( s ), 7 );
    [ 41 ]
    gap> ExtendTableHeadByPermCharValue( head, s, pi, pos );
    #I  after order 7 element from 2.B: have 191 classes
    gap> Last( head.SizesCentralizers ) = 7 * Size( CharacterTable( "He" ) );
    true
\end{verbatim}

By additional arguments, we find that the \texttt{7A} centralizer has the
structure $7 \times \He$,
and the normalizer has the structure $(7:3 \times \He).2$,
a subdirect product of $7:6$ and $\He.2$.

Since $\He$ has a pair of Galois conjugate classes of element order $17$,
we get also a pair of Galois conjugate classes of element order
$7 \cdot 17 = 119$.

\begin{verbatim}
    gap> ExtendTableHeadByCentralizerOrder( head, 119, 119, fail );
    #I  after order 119 element: have 192 classes
    gap> ExtendTableHeadByCentralizerOrder( head, 119, 119, fail );
    #I  after order 119 element: have 193 classes
\end{verbatim}

The second class of elements of order $7$, \texttt{7B},
is established by the fact that the subgroup $3.\Fi_{24}$ contains
two classes of elements of order $7$, with different values of the
degreee $196\,883$ character $\chi$ of $\M$.

\begin{verbatim}
    gap> u:= CharacterTable( "3.Fi24" );;
    gap> cand:= Filtered( Irr( u ), x -> x[1] <= 196883 );;
    gap> rest:= Sum( cand{ [ 1, 4, 5, 7, 8 ] } );;
    gap> pos:= Positions( OrdersClassRepresentatives( u ), 7 );
    [ 41, 43 ]
    gap> rest{ pos };
    [ 50, 1 ]
\end{verbatim}

Note that class $43$ fuses to \texttt{7B} because the restriction of $\chi$
to $2.\B$ has the value $50$ on \texttt{7A}.

\begin{verbatim}
    gap> ExtendTableHeadByCentralizerOrder( head, u, 7^5 * Factorial(7), [ 43 ] );
    #I  after order 7 element from 3.F3+.2: have 194 classes
\end{verbatim}

Now the sum of class lengths in \texttt{head} is equal to the order of $\M$.

\begin{verbatim}
    gap> Sum( head.SizesCentralizers, x -> head.Size / x ) = head.Size;
    true
\end{verbatim}

We initialize the character table head of $\M$.

\begin{verbatim}
    gap> m:= ConvertToCharacterTableNC( rec(
    >      UnderlyingCharacteristic:= 0,
    >      Size:= head.Size,
    >      SizesCentralizers:= head.SizesCentralizers,
    >      OrdersClassRepresentatives:= head.OrdersClassRepresentatives,
    >    ) );;
\end{verbatim}

\section{The power maps of $\M$}\label{Mpowermaps}

Using the element orders of the class representatives of the table head
of $\M$, and the partial class fusions from the subgroups used in the
previous sections, we compute approximations of the $p$-th power maps,
for primes $p$ up to the maximal element order in $\M$.

Note that we have not yet determined which of the two possible character
tables of the \texttt{3B} normalizer belongs to a subgroup of $\M$,
thus we exclude the corresponding partial fusion.

\begin{verbatim}
    gap> safe_fusions:= Filtered( head.fusions,
    >        r -> not IsIdenticalObj( r.subtable, facts[1] ) );;
    gap> Length( safe_fusions );
    6
\end{verbatim}

First we initialize the class fusions, compatible with the definitions of
the classes as given by the partial fusions which we have stored.

\begin{verbatim}
    gap> for r in safe_fusions do
    >      fus:= InitFusion( r.subtable, m );
    >      for i in [ 1 .. Length( r.map ) ] do
    >        if IsBound( r.map[i] ) then
    >          if IsInt( fus[i] ) then
    >            if fus[i] <> r.map[i] then
    >              Error( "fusion problem" );
    >            fi;
    >          elif IsInt( r.map[i] ) then
    >            if not r.map[i] in fus[i] then
    >              Error( "fusion problem" );
    >            fi;
    >          else
    >            if not IsSubset( fus[i], r.map[i] ) then
    >              Error( "fusion problem" );
    >            fi;
    >          fi;
    >          fus[i]:= r.map[i];
    >        fi;
    >      od;
    >      r.fus:= fus;
    >    od;
\end{verbatim}

Next we initialize approximations of the power maps of the table of $\M$,
and improve them using the compatibility of these maps with the power maps
of the subgroups w.~r.~t.~the current knowledge of the class fusions.
Note that also the knowledge about the class fusions increases this way.

\begin{verbatim}
    gap> maxorder:= Maximum( head.OrdersClassRepresentatives );
    119
    gap> powermaps:= [];;
    gap> primes:= Filtered( [ 1 .. maxorder ], IsPrimeInt );
    gap> for p in primes do
    >      powermaps[p]:= InitPowerMap( m, p );
    >      for r in safe_fusions do
    >        subpowermap:= PowerMap( r.subtable, p );
    >        if TransferDiagram( subpowermap, r.fus, powermaps[p] ) = fail then
    >          Error( "inconsistency" );
    >        fi;
    >      od;
    >    od;
\end{verbatim}

We repeat applying the compatibility conditions until no further
improvements are found.

\begin{verbatim}
    gap> found:= true;;
    gap> res:= "dummy";;  # avoid a syntax warning
    gap> while found do
    >      Print( "#I  start a round\n" );
    >      found:= false;
    >      for p in primes do
    >        for r in safe_fusions do
    >          subpowermap:= PowerMap( r.subtable, p );
    >          res:= TransferDiagram( subpowermap, r.fus, powermaps[p] );
    >          if res = fail then
    >            Error( "inconsistency" );
    >          elif ForAny( RecNames( res ), nam -> res.( nam ) <> [] ) then
    >            found:= true;
    >          fi;
    >        od;
    >      od;
    >    od;
    #I  start a round
    #I  start a round
    #I  start a round
    #I  start a round
\end{verbatim}

Let us see where the power maps are still not determined uniquely,
starting with the $5$-th power map.

\begin{verbatim}
    gap> pos:= PositionsProperty( powermaps[5], IsList );
    [ 157, 158, 163, 164, 187, 188, 189, 190, 192, 193 ]
    gap> head.OrdersClassRepresentatives{ pos };
    [ 15, 15, 39, 39, 59, 59, 71, 71, 119, 119 ]
\end{verbatim}

The ambiguities for the classes of the element orders $59$, $71$, and $119$
are understandable:
For each of these element orders, there is a pair of Galois conjugate
classes, and the subgroups whose class fusions we have used do not contain
these elements.

For each of the primes $l \in \{ 59, 71 \}$,
the field of $l$-th roots of unity contains a unique quadratic subfield,
which is $\Q(\sqrt{-l})$,
and the $p$-th power map, for $p$ coprime to $l$,
fixes a class of element order $l$ if and only if
the Galois automorphism that raises $l$-th roots of unity to the $p$-th power
fixes $\sqrt{-l}$.

In the case of element order $l = 119 = 7 \cdot 17$,
the field of $l$-th roots of unity contains the three quadratic subfields,
$\Q(\sqrt{-7})$, $\Q(\sqrt{17})$, and $\Q(\sqrt{-119})$.
In order to decide which of them actually occurs,
we look at a subgroup that contains elements of order $119$.
The \texttt{7A} centralizer in $\M$ has the structure $7 \times \He$,
and the normalizer has the structure $(7:3 \times \He).2$,
a subdirect product of $7:6$ and $\He.2$,
see Section~\ref{elements_7}.

The classes of element order $119$ in the normalizer correspond to
the classes of this element order in $\M$,
and the character values in the subgroup lie in the field $\Q(\sqrt{-119})$.

\begin{verbatim}
    gap> u:= CharacterTable( "(7:3xHe):2" );;
    gap> ConstructionInfoCharacterTable( u );
    [ "ConstructIndexTwoSubdirectProduct", "7:3", "7:6", "He", "He.2", 
      [ 117, 118, 119, 120, 121, 122, 123, 124, 125, 126, 127, 128, 129, 
          130, 131, 132, 133, 134, 135, 207, 208, 209, 210, 211, 212, 
          213, 214, 215, 216, 217, 218, 219, 220, 221, 222, 223, 224, 
          225, 297, 298, 299, 300, 301, 302, 303, 304, 305, 306, 307, 
          308, 309, 310, 311, 312, 313, 314, 315 ], (), () ]
    gap> pos:= Positions( OrdersClassRepresentatives( u ), 119 );
    [ 52, 53 ]
    gap> f:= Field( Rationals, List( Irr( u ), x -> x[pos[1]] ) );;
    gap> Sqrt(-119) in f;
    true
\end{verbatim}

We insert the relevant power map values.

\begin{verbatim}
    gap> for l in [ 59, 71, 119 ] do
    >      val:= Sqrt( -l );
    >      poss:= Positions( head.OrdersClassRepresentatives, l );
    >      for p in primes do
    >        if Gcd( l, p ) = 1 then
    >          if GaloisCyc( val, p ) = val then
    >            powermaps[p]{ poss }:= poss;
    >          else
    >            powermaps[p]{ poss }:= Reversed( poss );
    >          fi;
    >        fi;
    >      od;
    >    od;
\end{verbatim}

Now $p$-th power maps, for $p \geq 17$,
are determined uniquely except for the images of two classes
of element order $39$.
These classes had been found as roots of \texttt{3B} elements.

\begin{verbatim}
    gap> PositionsProperty( powermaps[17], IsList );
    [ 163, 164 ]
    gap> head.OrdersClassRepresentatives{ [ 163, 164 ] };
    [ 39, 39 ]
    gap> List( Filtered( head.fusions,
    >                    r -> IsSubset( r.map, [ 163, 164 ] ) ),
    >          r -> r.subtable );
    [ CharacterTable( "3^(1+12):6.Suz.2/[ 1, 19 ]" ) ]
\end{verbatim}

In order to decide whether the $p$-th power map fixes or swaps
the two classes, we consider elements of order $78$,
which are the square roots of the order $39$ elements.
There are three classes of element order $78$ in $\M$,
a rational class that powers to \texttt{2A}
and a pair of Galois conjugate classes that power to \texttt{2B}.

\begin{verbatim}
    gap> 78pos:= Positions( head.OrdersClassRepresentatives, 78 );
    [ 37, 132, 133 ]
    gap> head.fusions[1].subtable;
    CharacterTable( "2.B" )
    gap> Intersection( 78pos, head.fusions[1].map );
    [ 37 ]
    gap> s:= head.fusions[2].subtable;
    CharacterTable( "2^1+24.Co1" )
    gap> Intersection( 78pos, head.fusions[2].map );
    [ 132, 133 ]
    gap> Positions( head.fusions[2].map, 132 );
    [ 342 ]
    gap> Positions( head.fusions[2].map, 133 );
    [ 344 ]
    gap> PowerMap( s, 7 )[342];
    344
\end{verbatim}

Since the \texttt{3B} normalizer in $\M$ contains a pair of Galois conjugate
classes of element order $78$ which power to the generators of the normal
subgroup of order $3$,
these two classes fuse to the non-rational $\M$-classes of elements
of order $78$,
and their squares are the classes of element order $39$ we are interested in.

\begin{verbatim}
    gap> poss:= Filtered( head.fusions, r -> IsSubset( r.map, [ 163, 164 ] ) );;
    gap> List( poss, r -> r.subtable );
    [ CharacterTable( "3^(1+12):6.Suz.2/[ 1, 19 ]" ) ]
    gap> Position( poss[1].map, 163 );
    173
    gap> Position( poss[1].map, 164 );
    174
    gap> List( facts, s -> Positions( OrdersClassRepresentatives( s ), 39 ) );
    [ [ 173, 174 ], [ 173, 174 ] ]
    gap> List( facts, s -> PowerMap( s, 7 )[173] );
    [ 174, 174 ]
\end{verbatim}

The field of character values on the two classes of $\M$ is equal
to the corresponding field of character values in the \texttt{3B} normalizer,
which is $\Q(\sqrt{-39})$.
(Note that we have not yet decided which of the two candidate tables belong
to the \texttt{3B} normalizer,
but we get the same result for both candidates.)

\begin{verbatim}
    gap> fields:= List( facts,
    >                   s -> Field( Rationals, List( Irr( s ),
    >                                                x -> x[173] ) ) );;
    gap> Length( Set( fields ) );
    1
    gap> Sqrt(-39) in fields[1];
    true
\end{verbatim}

Now we can set the power map values on the two classes.

\begin{verbatim}
    gap> val:= Sqrt( -39 );;
    gap> poss:= [ 163, 164 ];;
    gap> for p in primes do
    >      if Gcd( 39, p ) = 1 then
    >        if GaloisCyc( val, p ) = val then
    >          powermaps[p]{ poss }:= poss;
    >        else
    >          powermaps[p]{ poss }:= Reversed( poss );
    >        fi;
    >      fi;
    >    od;
    gap> List( powermaps, Indeterminateness );
    [ , 2048, 1536,, 4,, 2,,,, 2,, 9,,,, 1,, 1,,,, 1,,,,,, 1,, 1,,,,,, 1,,
      ,, 1,, 1,,,, 1,,,,,, 1,,,,,, 1,, 1,,,,,, 1,,,, 1,, 1,,,,,, 1,,,, 1,,
      ,,,, 1,,,,,,,, 1,,,, 1,, 1,,,, 1,, 1,,,, 1 ]
\end{verbatim}

In the following,
we use the two candidates for the \texttt{3B} normalizer table
for answering most of the remaining questions about the power maps.
Again, the answers are equal for both candidate tables.

First we initialize the class fusion from the first candidate table \ldots

\begin{verbatim}
    gap> r:= First( head.fusions, r -> IsIdenticalObj( r.subtable, facts[1] ) );;
    gap> fus:= InitFusion( r.subtable, m );;
    gap> for i in [ 1 .. Length( r.map ) ] do
    >      if IsBound( r.map[i] ) then
    >        if IsInt( fus[i] ) then
    >          if fus[i] <> r.map[i] then
    >            Error( "fusion problem" );
    >          fi;
    >        elif IsInt( r.map[i] ) then
    >          if not r.map[i] in fus[i] then
    >            Error( "fusion problem" );
    >          fi;
    >        else
    >          if not IsSubset( fus[i], r.map[i] ) then
    >            Error( "fusion problem" );
    >          fi;
    >        fi;
    >        fus[i]:= r.map[i];
    >      fi;
    >    od;
    gap> r.fus:= fus;;
\end{verbatim}

\ldots and the class fusion from the second candidate table, \ldots

\begin{verbatim}
    gap> r2:= First( head2.fusions, r -> IsIdenticalObj( r.subtable, facts[2] ) );;
    gap> fus2:= InitFusion( r2.subtable, m );;
    gap> for i in [ 1 .. Length( r2.map ) ] do 
    >      if IsBound( r2.map[i] ) then
    >        if IsInt( fus2[i] ) then
    >          if fus2[i] <> r2.map[i] then
    >            Error( "fusion problem" );
    >          fi;
    >        elif IsInt( r2.map[i] ) then
    >          if not r2.map[i] in fus2[i] then
    >            Error( "fusion problem" );
    >          fi;
    >        else
    >          if not IsSubset( fus2[i], r2.map[i] ) then
    >            Error( "fusion problem" );
    >          fi;
    >        fi;
    >        fus2[i]:= r2.map[i];
    >      fi;
    >    od;
    gap> r2.fus:= fus2;;
\end{verbatim}

\ldots then we create an independent copy of the current approximations
of power maps, and apply the consistency conditions for class fusion and
power maps in the two cases.

\begin{verbatim}
    gap> powermaps2:= StructuralCopy( powermaps );;
    gap> s:= r.subtable;
    CharacterTable( "3^(1+12):6.Suz.2/[ 1, 19 ]" )
    gap> for p in primes do
    >      if TransferDiagram( PowerMap( s, p ), fus, powermaps[p] ) = fail then
    >        Error( "inconsistency" );
    >      fi;
    >    od;
    gap> s2:= r2.subtable;
    CharacterTable( "3^(1+12):6.Suz.2/[ 1, 20 ]" )
    gap> for p in primes do
    >      if TransferDiagram( PowerMap( s2, p ), fus2, powermaps2[p] ) = fail then
    >        Error( "inconsistency" );
    >      fi;
    >    od;
    gap> powermaps = powermaps2;
    true
    gap> List( powermaps, Indeterminateness );
    [ , 32, 64,, 1,, 1,,,, 1,, 1,,,, 1,, 1,,,, 1,,,,,, 1,, 1,,,,,, 1,,,, 
      1,, 1,,,, 1,,,,,, 1,,,,,, 1,, 1,,,,,, 1,,,, 1,, 1,,,,,, 1,,,, 1,,,,,
      , 1,,,,,,,, 1,,,, 1,, 1,,,, 1,, 1,,,, 1 ]
\end{verbatim}

One open question is about the squares of the non-rational classes
of element order $78$.

\begin{verbatim}
    gap> powermaps[2]{ [ 132, 133 ] };
    [ [ 163, 164 ], [ 163, 164 ] ]
    gap> pos78:= List( facts,
    >                  s -> Positions( OrdersClassRepresentatives( s ), 78 ) );
    [ [ 235, 236 ], [ 235, 236 ] ]
    gap> fus{ [ 235, 236 ] };
    [ [ 132, 133 ], [ 132, 133 ] ]
    gap> fus2{ [ 235, 236 ] };
    [ [ 132, 133 ], [ 132, 133 ] ]
\end{verbatim}

We may identify one class of element order $78$ in the \texttt{3B} normalizer
with the corresponding class of $\M$, and then draw conclusions.

\begin{verbatim}
    gap> fus[235]:= 132;;
    gap> fus2[235]:= 132;;
    gap> TransferDiagram( PowerMap( s, 2 ), fus, powermaps[2] ) <> fail;
    true
    gap> TransferDiagram( PowerMap( s2, 2 ), fus2, powermaps2[2] ) <> fail;
    true
    gap> powermaps = powermaps2;
    true
    gap> List( powermaps{ [ 2, 3 ] }, Indeterminateness );
    [ 8, 64 ]
\end{verbatim}

Since also the cubes of the concerned classes of element order $39$
are still not determined, this question is now decided using that the
2nd and the 3rd power map commute.

\begin{verbatim}
    gap> powermaps[3]{ [ 163, 164 ] };
    [ [ 183, 184 ], [ 183, 184 ] ]
    gap> TransferDiagram( powermaps[2], powermaps[3], powermaps[2] ) <> fail;
    true
    gap> List( powermaps{ [ 2, 3 ] }, Indeterminateness );
    [ 8, 16 ]
\end{verbatim}

The next open question is about the cubes of elements of order $93$.
The classes of element order $93$ --a pair of Galois conjugate classes--
have been found inside subgroups
of the type $S_3 \times \Th$, and they do not occur in other subgroups
we have considered.
Thus we may choose which of them cubes to the first class of element order
$31$.

\begin{verbatim}
    gap> poss:= Filtered( head.fusions,
    >                     r -> 93 in OrdersClassRepresentatives( r.subtable ) );;
    gap> List( poss, r -> r.subtable );
    [ CharacterTable( "ThxSym(3)" ) ]
    gap> pos93:= Positions( head.OrdersClassRepresentatives, 93 );
    [ 152, 153 ]
    gap> powermaps[3]{ pos93 };
    [ [ 179, 180 ], [ 179, 180 ] ]
    gap> powermaps[3][152]:= 179;;
    gap> TransferDiagram( PowerMap( poss[1].subtable, 3 ), poss[1].fus,
    >                     powermaps[3] ) <> fail;
    true
    gap> List( powermaps{ [ 2, 3 ] }, Indeterminateness );
    [ 8, 4 ]
\end{verbatim}

The next open question is about the cubes of elements of order $69$.
The classes of element order $69$ --a pair of Galois conjugate classes--
have been found inside subgroups
of the type $3.\Fi_{24}$, and they do not occur in other subgroups
we have considered.
Thus we may choose which of them cubes to the first class of element order
$23$.

\begin{verbatim}
    gap> poss:= Filtered( head.fusions, 
    >                     r -> 69 in OrdersClassRepresentatives( r.subtable ) );;
    gap> List( poss, r -> r.subtable );
    [ CharacterTable( "3.F3+.2" ) ]
    gap> pos69:= Positions( head.OrdersClassRepresentatives, 69 );
    [ 142, 143 ]
    gap> powermaps[3]{ pos69 };
    [ [ 177, 178 ], [ 177, 178 ] ]
    gap> powermaps[3]{ [ 142, 143 ] }:= [ 177, 178 ];;
    gap> TransferDiagram( PowerMap( poss[1].subtable, 3 ), poss[1].fus,
    >                     powermaps[3] ) <> fail;
    true
    gap> List( powermaps{ [ 2, 3 ] }, Indeterminateness );
    [ 8, 1 ]
\end{verbatim}

The next open question is about the squares of certain elements of order $46$.
There are two pairs of Galois conjugate classes of element order $46$,
and the 2nd power map is not yet determined for those classes which power
to the class \texttt{2B}.

\begin{verbatim}
    gap> pos46:= Positions( head.OrdersClassRepresentatives, 46 );
    [ 26, 27, 118, 120 ]
    gap> powermaps[2]{ pos46 };
    [ 177, 178, [ 177, 178 ], [ 177, 178 ] ]
    gap> powermaps[23]{ pos46 };
    [ 2, 2, 44, 44 ]
\end{verbatim}

We have defined the two classes of element order $23$ as squares of those
two classes of element order $46$ that power to \texttt{2A},
and we have not yet distinguished the other two classes of element order $46$.
Thus we may set the power map values.

\begin{verbatim}
    gap> powermaps[2]{ [ 118, 120 ] }:= [ 177, 178 ];;
    gap> Indeterminateness ( powermaps[2] );
    2
\end{verbatim}

Now just one value is left to be determined,
the square of a class of element order $18$ and centralizer order $3888$.

\begin{verbatim}
    gap> powermaps[2][78];
    [ 155, 156 ]
    gap> head.OrdersClassRepresentatives[78];
    18
    gap> head.SizesCentralizers[78];
    3888
\end{verbatim}

There are two classes with this property in $\M$,
both are roots of the generators of the normal subgroup of order $3$ in
$3^{1+12}_+.2.\Suz.2$, and the corresponding two classes in this subgroup
have the same square.

\begin{verbatim}
    gap> Filtered( [ 1 .. Length( head.OrdersClassRepresentatives ) ],
    >              i -> head.OrdersClassRepresentatives[i] = 18 and
    >                   head.SizesCentralizers[i] = 3888 );
    [ 78, 79 ]
    gap> powermaps[3]{ [ 78, 79 ] };
    [ 52, 52 ]
    gap> powermaps[2][52];
    154
    gap> First( head.fusions, r -> 154 in r.map ).subtable;
    CharacterTable( "3^(1+12):6.Suz.2/[ 1, 19 ]" )
    gap> s:= facts[1];
    CharacterTable( "3^(1+12):6.Suz.2/[ 1, 19 ]" )
    gap> pos18:= Filtered( [ 1 .. NrConjugacyClasses( s ) ],
    >                i -> OrdersClassRepresentatives( s )[i] = 18 and
    >                     SizesCentralizers( s )[i] = 3888 );
    [ 67, 83 ]
    gap> PowerMap( s, 2 ){ pos18 };
    [ 24, 24 ]
    gap> s:= facts[2];
    CharacterTable( "3^(1+12):6.Suz.2/[ 1, 20 ]" )
    gap> pos18:= Filtered( [ 1 .. NrConjugacyClasses( s ) ],
    >                i -> OrdersClassRepresentatives( s )[i] = 18 and
    >                     SizesCentralizers( s )[i] = 3888 );
    [ 67, 83 ]
    gap> PowerMap( s, 2 ){ pos18 };
    [ 24, 24 ]
\end{verbatim}

We set the last missing value,
and improve the approximations of the class fusions we have used,
by applying the consistency criteria.

\begin{verbatim}
    gap> powermaps[2][78]:= powermaps[2][79];;
    gap> for r in safe_fusions do
    >      if not TestConsistencyMaps( ComputedPowerMaps( r.subtable ), r.fus,
    >                                  powermaps ) then
    >        Error( "inconsistent!" );
    >      fi;
    >    od;
    gap> r:= First( head.fusions, r -> IsIdenticalObj( r.subtable, facts[1] ) );;
    gap> TestConsistencyMaps( ComputedPowerMaps( r.subtable ), r.fus,
    >                         powermaps );
    true
    gap> r2:= First( head2.fusions, r -> IsIdenticalObj( r.subtable, facts[2] ) );;
    gap> TestConsistencyMaps( ComputedPowerMaps( r2.subtable ), r2.fus,
    >                         powermaps );
    true
    gap> SetComputedPowerMaps( m, powermaps );
\end{verbatim}

\section{The degree $196\,883$ character $\chi$ of $\M$}\label{sect:natcharM}

We know the values of the irreducible degree $196\,883$ character $\chi$
of $\M$ on the classes of $2.\B$, by Section~\ref{natural}.

\begin{verbatim}
    gap> r:= head.fusions[1];;
    gap> s:= r.subtable;
    CharacterTable( "2.B" )
    gap> cand:= Filtered( Irr( s ), x -> x[1] <= 196883 );;
    gap> List( cand, x -> x[1] );
    [ 1, 4371, 96255, 96256 ]
    gap> rest:= Sum( cand );;
    gap> rest[1];
    196883
\end{verbatim}

Thus we know the values of $\chi$ on those classes of $\M$
that are known as images of the class fusion from $2.\B$.
This yields $111$ out of the $194$ character values.

\begin{verbatim}
    gap> chi:= [];;
    gap> map:= r.fus;;
    gap> for i in [ 1 .. Length( map ) ] do
    >      if IsInt( map[i] ) then
    >        chi[ map[i] ]:= rest[i];
    >      fi;
    >    od;
    gap> Number( chi );
    111
\end{verbatim}

Also the restriction of $\chi$ to $3.\Fi_{24}$ is known,
by Section~\ref{natural}.
This yields $29$ more character values.

\begin{verbatim}
    gap> r:= head.fusions[3];;
    gap> s:= r.subtable;
    CharacterTable( "3.F3+.2" )
    gap> cand:= Filtered( Irr( s ), x -> x[1] <= 196883 );;
    gap> rest:= Sum( cand{ [ 1, 4, 5, 7, 8 ] } );;
    gap> rest[1];
    196883
    gap> map:= r.fus;;
    gap> for i in [ 1 .. Length( map ) ] do
    >      if IsInt( map[i] ) then
    >        if IsBound( chi[ map[i] ] ) and chi[ map[i] ] <> rest[i] then
    >          Error( "inconsistency!" );
    >        fi;
    >        chi[ map[i] ]:= rest[i];
    >      fi;
    >    od;
    gap> Number( chi );
    140
\end{verbatim}

Now we compute the restriction of $\chi$ to the \texttt{2B} normalizer.
There are only $13$ possible irreducible constituents of this restriction.
We consider the matrix of values of these characters on those classes
for which the class fusion to $\M$ is uniquely known \emph{and}
the value of $\chi$ on the image class is known.
This matrix has full rank, thus we can directly compute the decomposition
of the restriction into irreducibles,
and get $41$ more character values.

\begin{verbatim}
    gap> r:= head.fusions[2];;
    gap> s:= r.subtable;
    CharacterTable( "2^1+24.Co1" )
    gap> cand:= Filtered( Irr( s ), x -> x[1] <= chi[1] );;
    gap> map:= r.fus;;
    gap> knownpos:= Filtered( [ 1 .. Length( map ) ],
    >                      i -> IsInt( map[i] ) and IsBound( chi[ map[i] ] ) );;
    gap> rest:= List( knownpos, i -> chi[ map[i] ] );;
    gap> mat:= List( cand, x -> x{ knownpos } );;
    gap> Length( mat );
    13
    gap> RankMat( mat );
    13
    gap> sol:= SolutionMat( mat, rest );
    [ 0, 0, 1, 0, 0, 0, 0, 0, 0, 0, 0, 1, 1 ]
    gap> rest:= sol * cand;;
    gap> for i in [ 1 .. Length( map ) ] do
    >      if IsInt( map[i] ) then chi[ map[i] ]:= rest[i]; fi;
    >    od;
    gap> Number( chi );
    181
\end{verbatim}

Which values are still missing?

\begin{verbatim}
    gap> missing:= Filtered( [ 1..194 ], i -> not IsBound( chi[i] ) );
    [ 151, 152, 153, 160, 169, 170, 186, 187, 188, 189, 190, 192, 193 ]
    gap> head.OrdersClassRepresentatives{ missing };
    [ 57, 93, 93, 27, 95, 95, 41, 59, 59, 71, 71, 119, 119 ]
    gap> head.SizesCentralizers{ missing };
    [ 57, 93, 93, 243, 95, 95, 41, 59, 59, 71, 71, 119, 119 ]
\end{verbatim}

For $g \in \M$, we have $\chi(g) \equiv \chi(g^p) \pmod{p}$
and $|\chi(g)|^2 < |C_{\M}(g)|$.
We apply these conditions.
In all cases except one, the centralizer orders are small enough
for determining the character value uniquely.

\begin{verbatim}
    gap> for i in missing do
    >      ord:= head.OrdersClassRepresentatives[i];
    >      divs:= PrimeDivisors( ord );
    >      if ForAll( divs, p -> IsBound( chi[ powermaps[p][i] ] ) ) then
    >        congr:= List( divs, p -> chi[ powermaps[p][i] ] mod p );
    >        res:= ChineseRem( divs, congr );
    >        modulus:= Lcm( divs );
    >        c:= head.SizesCentralizers[i];
    >        Print( "#I  |g| = ", head.OrdersClassRepresentatives[i],
    >               ", |C_M(g)| = ", c,
    >               ": value ", res, " modulo ", modulus, "\n" );
    >        if ( res + 2 * modulus )^2 >= c and ( res - 2 * modulus )^2 >= c then
    >          cand:= Filtered( res + [ -1 .. 1 ] * modulus, a -> a^2 < c );
    >          if Length( cand ) = 1 then
    >            chi[i]:= cand[1];
    >          fi;
    >        fi;
    >      fi;
    >    od;
    #I  |g| = 57, |C_M(g)| = 57: value 56 modulo 57
    #I  |g| = 93, |C_M(g)| = 93: value 92 modulo 93
    #I  |g| = 93, |C_M(g)| = 93: value 92 modulo 93
    #I  |g| = 27, |C_M(g)| = 243: value 2 modulo 3
    #I  |g| = 95, |C_M(g)| = 95: value 0 modulo 95
    #I  |g| = 95, |C_M(g)| = 95: value 0 modulo 95
    #I  |g| = 41, |C_M(g)| = 41: value 1 modulo 41
    #I  |g| = 59, |C_M(g)| = 59: value 0 modulo 59
    #I  |g| = 59, |C_M(g)| = 59: value 0 modulo 59
    #I  |g| = 71, |C_M(g)| = 71: value 0 modulo 71
    #I  |g| = 71, |C_M(g)| = 71: value 0 modulo 71
    #I  |g| = 119, |C_M(g)| = 119: value 118 modulo 119
    #I  |g| = 119, |C_M(g)| = 119: value 118 modulo 119
    gap> missing:= Filtered( [ 1..194 ], i -> not IsBound( chi[i] ) );
    [ 160 ]
\end{verbatim}

The one missing value can be computed from the scalar product
with the trivial character.

\begin{verbatim}
    gap> diff:= Difference( [ 1 .. NrConjugacyClasses( m ) ], missing );;
    gap> classes:= SizesConjugacyClasses( m );;
    gap> sum:= Sum( diff, i -> classes[i] * chi[i] );
    -6650349175263480459970863415322722279882752000000000
    gap> chi[ missing[1] ]:= - sum / classes[ missing[1] ];
    2
\end{verbatim}

Now we decide which of the two candidates for the character table of
the \texttt{3B} normalizer is the correct one.
For the first candidate, the restriction of $\chi$
cannot be decomposed into irreducibles.

\begin{verbatim}
    gap> r:= First( head.fusions, r -> IsIdenticalObj( r.subtable, facts[1] ) );;
    gap> map:= r.fus;;
    gap> knownpos:= Filtered( [ 1 .. Length( map ) ], i -> IsInt( map[i] ) );;
    gap> rest:= List( knownpos, i -> chi[ map[i] ] );;
    gap> cand:= Filtered( Irr( r.subtable ), x -> x[1] <= chi[1] );;
    gap> mat:= List( cand, x -> x{ knownpos } );;
    gap> Length( mat );
    95
    gap> RankMat( mat );
    88
    gap> SolutionMat( mat, rest );
    fail
\end{verbatim}

The second candidate admits a decomposition.

\begin{verbatim}
    gap> r2:= First( head2.fusions, r -> IsIdenticalObj( r.subtable, facts[2] ) );;
    gap> map:= r2.fus;;
    gap> knownpos:= Filtered( [ 1 .. Length( map ) ], i -> IsInt( map[i] ) );;
    gap> rest:= List( knownpos, i -> chi[ map[i] ] );;
    gap> cand:= Filtered( Irr( r2.subtable ), x -> x[1] <= chi[1] );;
    gap> mat:= List( cand, x -> x{ knownpos } );;
    gap> Length( mat );
    95
    gap> RankMat( mat );
    88
    gap> SolutionMat( mat, rest );
    [ 0, 0, 0, 1, 0, 0, 0, 0, 0, 0, 0, 0, 0, 0, 0, 0, 0, 0, 0, 0, 0, 0, 
      0, 0, 0, 0, 0, 0, 0, 0, 0, 0, 0, 0, 0, 0, 0, 0, 0, 0, 0, 0, 0, 0, 
      0, 0, 0, 0, 0, 0, 0, 0, 0, 0, 0, 0, 0, 0, 0, 0, 0, 0, 0, 0, 0, 0, 
      0, 0, 0, 0, 0, 0, 0, 0, 0, 0, 0, 0, 0, 0, 0, 0, 0, 0, 0, 0, 0, 0, 
      0, 0, 0, 1, 1, 1, 0 ]
    gap> Add( safe_fusions, r2 );
\end{verbatim}

The character table of the second candidate is equivalent to
the character table that is stored in the {\GAP} character table library.

\begin{verbatim}
    gap> TransformingPermutationsCharacterTables( r2.subtable,
    >        CharacterTable( "MN3B" ) ) <> fail;
    true
\end{verbatim}

\section{The irreducible characters of $\M$}\label{sect:irreduciblesM}

We will not compute the irreducibles of $\M$ from scratch but verify the
irreducibles from the {\ATLAS} character table of $\M$,
in the sense that we use the characters printed in the {\ATLAS} as an
``oracle''.
For that,
we compute first a bijection between the columns of our character table head
and those of the {\ATLAS} character table of $\M$.
This is done by using the following invariants:
element orders, centralizer orders,
the values of $\chi$, and the indirection of $\chi$ by the $2$nd power map.

\begin{verbatim}
    gap> invs:= TransposedMat( [
    >      OrdersClassRepresentatives( m ),
    >      SizesCentralizers( m ),
    >      chi,
    >      CompositionMaps( chi, PowerMap( m, 2 ) ) ] );;
    gap> invs_set:= Set( invs );;
    gap> Length( invs_set );
    172
    gap> atlas_m:= CharacterTable( "M" );;
    gap> invs_atlas:= TransposedMat( [
    >      OrdersClassRepresentatives( atlas_m ),
    >      SizesCentralizers( atlas_m ),
    >      Irr( atlas_m )[2],
    >      CompositionMaps( Irr( atlas_m )[2], PowerMap( atlas_m, 2 ) ) ] );;
    gap> invs_atlas_set:= Set( invs_atlas );;
    gap> invs_atlas_set = invs_set;
    true
\end{verbatim}

In particular, we see that the sets of invariants are equal for the two
tables.

Note that we cannot get a better choice of invariants,
since there are $22$ pairs of Galois conjugate classes in $\M$,
and our current knowledge does not allow us to distinguish
the classes of each pair.

Now we compute a permutation that maps the classes of the {\ATLAS} table
to suitable classes of our table head,
permute the irreducibles of the {\ATLAS} table accordingly,
and create the ``oracle'' list.

(The explicit permutation $(32,33)(179,180)$ makes sure that the power maps
of the {\ATLAS} table and of our table are compatible.)

\begin{verbatim}
    gap> pi1:= SortingPerm( invs );;
    gap> pi2:= SortingPerm( invs_atlas );;
    gap> pi:= pi2 / pi1 * (32,33)(179,180);;
    gap> oracle:= List( Irr( atlas_m ), x -> Permuted( x, pi ) );;
\end{verbatim}

In order to prohibit that {\GAP} tries to compute table automorphisms
of our table head of $\M$ (which is impossible without knowing the
irreducible characters),
we set a trivial group as value of the attribute
\texttt{AutomorphismsOfTable};
this will be revised as soon as the irreducibles are known.

\begin{verbatim}
    gap> SetAutomorphismsOfTable( m, Group( () ) );
\end{verbatim}


First we compute candidates for the class fusion from $2.\B$,
starting from the approximation we have already computed.
Before we apply {\GAP}'s criteria for computing possible class fusions,
we decide about the images of four classes of element orders $40$ and $44$.

\begin{verbatim}
    gap> r:= safe_fusions[1];;
    gap> s:= r.subtable;
    CharacterTable( "2.B" )
    gap> pos:= [ 217, 218, 222, 223 ];;
    gap> r.fus{ pos };
    [ [ 110, 111 ], [ 110, 111 ], [ 88, 89 ], [ 88, 89 ] ]
    gap> OrdersClassRepresentatives( s ){ pos };
    [ 40, 40, 44, 44 ]
\end{verbatim}

The two classes of element order $40$ are a pair of Galois conjugates
in both $2.\B$ and $\M$.
Note that their elements are roots of \texttt{2B} elements in $\M$,
and the classes $110$ and $111$ correspond to non-rational elements
of order $40$ in the \texttt{2B} normalizer.

\begin{verbatim}
    gap> r.fus[ PowerMap( s, 20 )[217] ];
    44
    gap> r2:= safe_fusions[2];;
    gap> s2:= r2.subtable;
    CharacterTable( "2^1+24.Co1" )
    gap> Position( r2.map, 110 );
    273
    gap> ForAll( Irr( s2 ), x -> IsInt( x[273] ) );
    false
\end{verbatim}

We may freely choose the fusion from the classes $217$ and $218$ of $2.\B$
because there is a table automorphism of $2.\B$ that swaps exactly these
two classes.

\begin{verbatim}
    gap> (217,218) in AutomorphismsOfTable( s );
    true
    gap> r.fus{ [ 217, 218 ] }:= [ 110, 111 ];;
\end{verbatim}

With the same argument,
also the two classes of element order $44$ are a pair of Galois conjugates
both in $2.\B$ and $\M$. 

\begin{verbatim}
    gap> r.fus[ PowerMap( s, 22 )[ 222 ] ];
    44
    gap> Position( r2.map, 88 );
    178
    gap> ForAll( Irr( s2 ), x -> IsInt( x[178] ) );
    false
\end{verbatim}

There is a table automorphism of $2.\B$ that swaps three pairs of classes,
where the classes of element order $44$ form one pair,
and each of the other two pairs is fused in $\M$.
Thus we may again freely choose the fusion from $222$ and $223$.

\begin{verbatim}
    gap> (143,144)(222,223)(244,245) in AutomorphismsOfTable( s );
    true
    gap> r.fus{ [ 143, 144, 244, 245 ] };
    [ 15, 15, 34, 34 ]
    gap> r.fus{ [ 222, 223 ] }:= [ 88, 89 ];;
\end{verbatim}

Now the remaining open questions about the class fusion from $2.\B$
can be answered by {\GAP}'s function \texttt{PossibleClassFusions}.

\begin{verbatim}
    gap> knownirr:= [ TrivialCharacter( m ), chi ];;
    gap> poss:= PossibleClassFusions( s, m,
    >               rec( chars:= knownirr, fusionmap:= r.fus ) );;
    gap> List( poss, Indeterminateness );
    [ 1 ]
\end{verbatim}

Now we can induce the irreducibles of $2.\B$ to $\M$.

\begin{verbatim}
    gap> induced:= InducedClassFunctionsByFusionMap( s, m, Irr( s ), poss[1] );;
\end{verbatim}

Next we compute candidates for the class fusions from the subgroups
$2^{1+24}_+.\Co_1$, $3^{1+12}_+.2.\Suz.2$, and $3.\Fi_{24}$.
Here we enter also the characters of $\M$ obtained by induction from $2.\B$,
because their restrictions to the subgroups provide additional conditions.

The fusion from $2^{1+24}_+.\Co_1$ is determined uniquely up to automorphisms
of the subgroup table.
We extend the list of known induced characters.

\begin{verbatim}
    gap> poss:= PossibleClassFusions( s2, m,
    >               rec( chars:= Concatenation( knownirr, induced ),
    >                    fusionmap:= r2.fus ) );;
    gap> List( poss, Indeterminateness );
    [ 1 ]
    gap> Append( induced,
    >        InducedClassFunctionsByFusionMap( s2, m, Irr( s2 ), poss[1] ) );
\end{verbatim}

In order to compute the fusion from $3^{1+12}_+.2.\Suz.2$,
we have to consider two classes of element order $56$ first.
They are a pair of Galois conjugates both in $3^{1+12}_+.2.\Suz.2$ and $\M$,
and we may freely choose their fusion because there is a table automorphism
of $3^{1+12}_+.2.\Suz.2$ that swaps exactly these two classes.

\begin{verbatim}
    gap> r:= safe_fusions[7];;
    gap> s:= r.subtable;
    CharacterTable( "3^(1+12):6.Suz.2/[ 1, 20 ]" )
    gap> pos:= Positions( OrdersClassRepresentatives( s ), 56 );
    [ 250, 251 ]
    gap> r.fus{ pos };
    [ [ 125, 126 ], [ 125, 126 ] ]
    gap> r.fus[ PowerMap( s, 28 )[ 250 ] ];
    44
    gap> Position( r2.map, 125 );
    319
    gap> ForAll( Irr( s2 ), x -> IsInt( x[319] ) );
    false
    gap> (250, 251) in AutomorphismsOfTable( s );
    true
    gap> r.fus{ [ 250, 251 ] }:= [ 125, 126 ];;
\end{verbatim}

Now the class fusion to $\M$ is determined uniquely,
and we extend the list of induced characters.

\begin{verbatim}
    gap> poss:= PossibleClassFusions( s, m,
    >               rec( chars:= Concatenation( knownirr, induced ),
    >                    fusionmap:= r.fus ) );;
    gap> List( poss, Indeterminateness );
    [ 1 ]
    gap> Append( induced,
    >        InducedClassFunctionsByFusionMap( s, m, Irr( s ), poss[1] ) );
\end{verbatim}

The fusion from $3.\Fi_{24}$ is determined uniquely up to automorphisms
of the subgroup table.
We extend the list of known induced characters.

\begin{verbatim}
    gap> r:= safe_fusions[3];;
    gap> s:= r.subtable;
    CharacterTable( "3.F3+.2" )
    gap> poss:= PossibleClassFusions( s, m,
    >               rec( chars:= Concatenation( knownirr, induced ),
    >                    fusionmap:= r.fus ) );;
    gap> List( poss, Indeterminateness );
    [ 1 ]
    gap> Append( induced,
    >        InducedClassFunctionsByFusionMap( s, m, Irr( s ), poss[1] ) );
\end{verbatim}

Next we induce the irreducible characters of cyclic subgroups.

\begin{verbatim}
    gap> Append( induced,
    >      InducedCyclic( m, [ 2 .. NrConjugacyClasses( m ) ], "all" ) );
\end{verbatim}

Now we reduce the induced characters with the two known irreducibles of $\M$,
and apply the LLL algorithm to the result of the reduction;
this yields four new irreducibles.

\begin{verbatim}
    gap> red:= Reduced( m, knownirr, induced );;
    gap> Length( red.irreducibles );
    0
    gap> lll:= LLL( m, red.remainders );;
    gap> Length( lll.irreducibles );
    4
\end{verbatim}

We extend the list of known irreducibles,
reduce the induced characters,
and apply LLL again.

\begin{verbatim}
    gap> knownirr:= Union( knownirr, lll.irreducibles );;
    gap> red:= Reduced( m, knownirr, induced );;
    gap> Length( red.irreducibles );
    0
    gap> lll:= LLL( m, red.remainders );;
    gap> Length( lll.irreducibles );
    0
\end{verbatim}

Now we use the irreducibles of the {\ATLAS} table of $\M$ as an oracle,
as follows.
Whenever a character from the oracle list belongs to the $\Z$-lattice
that is spanned by \texttt{lll.remainders}
then we regard this character as verified,
since we can compute the coefficients of the $\Z$-linear combination,
form the character, and check that it has indeed norm $1$.

\begin{verbatim}
    gap> mat:= MatScalarProducts( m, oracle, lll.remainders );;
    gap> norm:= NormalFormIntMat( mat, 4 );;
    gap> rowtrans:= norm.rowtrans;;
    gap> normal:= norm.normal{ [ 1 .. norm.rank ] };;
    gap> one:= IdentityMat( NrConjugacyClasses( m ) );;
    gap> for i in [ 2 .. Length( one ) ] do
    >      extmat:= Concatenation( normal, [ one[i] ] );
    >      extlen:= Length( extmat );
    >      extnorm:= NormalFormIntMat( extmat, 4 );
    >      if extnorm.rank = Length( extnorm.normal ) or
    >         extnorm.rowtrans[ extlen ][ extlen ] <> 1 then
    >        coeffs:= fail;
    >      else
    >        coeffs:= - extnorm.rowtrans[ extlen ]{ [ 1 .. extnorm.rank ] }
    >                   * rowtrans{ [ 1 .. extnorm.rank ] };
    >      fi;
    >      if coeffs <> fail and ForAll( coeffs, IsInt ) then
    >        # The vector lies in the lattice.
    >        chi:= coeffs * lll.remainders;
    >        if not chi in knownirr then
    >          Add( knownirr, chi );
    >        fi;
    >      fi;
    >    od;
    gap> Length( knownirr );
    66
    gap> Set( knownirr, chi -> ScalarProduct( m, chi, chi ) );
    [ 1 ]
\end{verbatim}

We take the generators of the $\Z$-lattice and some symmetrizations
of the known irreducibles,
reuce them with the known irreducibles,
and apply LLL again.

\begin{verbatim}
    gap> red:= Reduced( m, knownirr, lll.remainders );;
    gap> Length( red.irreducibles );
    0
    gap> sym:= Symmetrizations( m, knownirr, 2 );;
    gap> sym:= Reduced( m, knownirr, sym );;
    gap> Length( sym.irreducibles );
    0
    gap> lll:= LLL( m, Concatenation( red.remainders, sym.remainders ) );;
    gap> Length( lll.irreducibles );
    0
\end{verbatim}

We use the above oracle again, for the new $\Z$-lattice.

\begin{verbatim}
    gap> mat:= MatScalarProducts( m, oracle, lll.remainders );;
    gap> norm:= NormalFormIntMat( mat, 4 );;
    gap> rowtrans:= norm.rowtrans;;
    gap> normal:= norm.normal{ [ 1 .. norm.rank ] };;
    gap> one:= IdentityMat( NrConjugacyClasses( m ) );;
    gap> for i in [ 2 .. Length( one ) ] do
    >      extmat:= Concatenation( normal, [ one[i] ] );
    >      extlen:= Length( extmat );
    >      extnorm:= NormalFormIntMat( extmat, 4 );
    >      if extnorm.rank = Length( extnorm.normal ) or
    >         extnorm.rowtrans[ extlen ][ extlen ] <> 1 then
    >        coeffs:= fail;
    >      else
    >        coeffs:= - extnorm.rowtrans[ extlen ]{ [ 1 .. extnorm.rank ] }
    >                   * rowtrans{ [ 1 .. extnorm.rank ] };
    >      fi;
    >      if coeffs <> fail and ForAll( coeffs, IsInt ) then
    >        Add( knownirr, coeffs * lll.remainders );
    >      fi;
    >    od;
    gap> Length( knownirr );
    194
\end{verbatim}

Now we are done.
As stated in the beginning of this section,
we unbind the stored trivial value for \texttt{AutomorphismsOfTable}.

\begin{verbatim}
    gap> SetIrr( m, List( knownirr, x -> ClassFunction( m, x ) ) );
    gap> ResetFilterObj( m, HasAutomorphismsOfTable );
    gap> TransformingPermutationsCharacterTables( m, atlas_m ) <> fail;
    true
\end{verbatim}

\section{Appendix: The character table of $2^{1+24}_+.\Co_1$}%
\label{sect:table_c2b}

The centralizer $C$ of a \texttt{2B} element in ${\M}$
has the structure $2^{1+24}_+.\Co_1$,
which can be constructed as follows.

Consider a subdirect product $H$ of two groups $H/X$
and $H/N$,
where $X$ is a cyclic group of order two,
$N$ is an extraspecial group $2^{1+24}_+$,
$H/N$ is the double cover of $\Co_1$, and
$H/X$ is an extension of $2^{1+24}_+$ by $\Co_1$.

The centre of $H$ is a Klein four group $E$ whose order two subgroups
are $X$, $Y = Z(N)$, and a third subgroup $D$.
We have $C \cong H/D$.

%

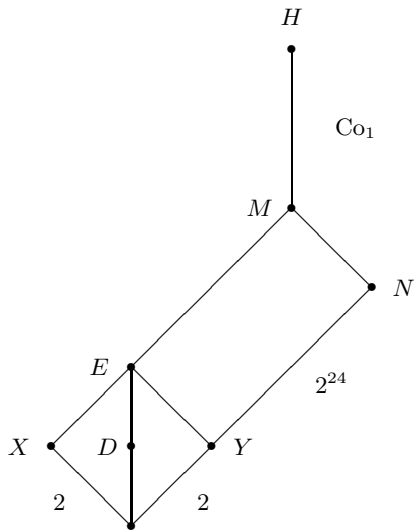
\begin{center}
\tthdump{\setlength{\unitlength}{3pt}}
\begin{picture}(50,70)
\put(20, 0){\circle*{1}}
\put(10,10){\circle*{1}} \put(6,10){\makebox(0,0){$X$}}
\put(20,10){\circle*{1}} \put(17,10){\makebox(0,0){$D$}}
\put(30,10){\circle*{1}} \put(34,10){\makebox(0,0){$Y$}}
\put(20,20){\circle*{1}} \put(16,20){\makebox(0,0){$E$}}
\put(20, 0){\line(-1,1){10}}
\put(20, 0){\line(0,1){20}}
\put(20, 0){\line(1,1){10}}
\put(20,20){\line(-1,-1){10}}
\put(20,20){\line(1,-1){10}}
\put(40,40){\circle*{1}} \put(36,40){\makebox(0,0){$M$}}
\put(50,30){\circle*{1}} \put(54,30){\makebox(0,0){$N$}}
\put(40,60){\circle*{1}} \put(40,64){\makebox(0,0){$H$}}
\put(20,20){\line(1,1){20}}
\put(30,10){\line(1,1){20}}
\put(50,30){\line(-1,1){10}}
\put(40,40){\line(0,1){20}}
\put(11,3){\makebox(0,0){$2$}}
\put(29,3){\makebox(0,0){$2$}}
\put(45,18){\makebox(0,0){$2^{24}$}}
\put(48,50){\makebox(0,0){$\Co_1$}}
\end{picture}
\end{center}

{}From the $2$-local construction of a matrix representation for ${\M}$,
we know the following faithful representations,
given by generators that are preimages of standard generators of $\Co_1$.

\begin{itemize}
\item
   A monomial permutation representation of $H/E \cong 2^{24}.\Co_1$,
   of degree $98\,280$ over the field with three elements,
\item
   a matrix representation of $H/X \cong 2^{1+24}_+.\Co_1$,
   of dimension $4\,096$ over the field with three elements,
\item
   a matrix representation of $H/N \cong 2.\Co_1$,
   of dimension $24$ over the field with three elements.
\end{itemize}


We proceed in the following steps.

\begin{itemize}
\item
   First we compute conjugacy class representatives of $H/E$.
\item
   A faithful permutation representation of $H/Y$ on
   $2 \cdot 196\,560 = 393\,120$ points is obtained by
   glueing the permutation generators of $H/E$
   (on $2 \cdot 98\,280 = 196\,560$ points)
   and $H/N$
   (the smallest permutation representation of $2.\Co_1$,
   on $196\,560$ points) together.

   The character table of this permutation group can be computed
   directly with {\MAGMA} in about $16$ hours of CPU time.
\item
   Starting from the character table of the factor group $H/E$,
   the $3$-modular matrix representation of dimension $4\,096$ of $H/X$ is used
   to compute necessary class splittings from $H/E$ to $H/X$, as follows.

   This representation lifts to characteristic zero because its restriction
   to the extraspecial group $M/X$
   is the unique faithful irreducible $3$-modular representation of $M/X$,
   and because this representation extends to the full automorphism group
   of the extraspecial group.
   Thus we can compute the Brauer character values of the representation
   on the $3$-regular classes of $H/X$,
   and interpret the values as those of an ordinary character $\psi$, say.
   The tensor square $\psi^2$ belongs to the group $H/E$,
   and the known values of $\psi$ suffice to determine the decomposition
   of $\psi^2$ into irreducibles of $H/E$,
   and thus to compute also the values of $\psi^2$ on $3$-singular classes.
   Taking square roots, we get all values of $\psi$, up to signs.
   (We cannot distinguish which of the two values belongs to which of the two
   preimage classes, which just means that we are defining these classes
   by choosing the positive value for one of them, and the negative value
   for the other one.)
\item
   From now on, we argue character-theoretically.

   The missing irreducible characters of $H/X$ are computed as
   tensor products of $\psi$ with the irreducible characters of
   the factor group $H/M \cong \Co_1$.

   Note that this procedure yields enough new irreducible characters
   such that the sum of degree squares of all known irreducibles of $H/X$
   equals the order of this group.
   This implies that there are not more class splittings w.~r.~t.~the
   fusion from $H/X$ to $H/E$ than the splittings forced by the values
   of $\psi$.
\item
   Using the two class splittings from $H/X$ and $H/Y$ to $H/E$,
   we compute necessary class splittings from $H$ to $H/E$.
   That is, we create a character table head for $H$ together with
   class fusions to $H/X$ and $H/Y$,
   assuming that not more columns occur than is forced by $H/X$ and $H/Y$:
   For each class of $H/E$ that splits in both $H/X$ and $H/Y$,
   we get four preimage classes in $H$.
   For each class that splits in exactly one of $H/X$ and $H/Y$,
   we get two preimage classes in $H$.
   For each class that splits in none of $H/X$ and $H/Y$,
   we get one preimage class in $H$.

   Then we take those irreducible characters of $H/X$ and $H/Y$,
   respectively, that do not have $E/X$ or $E/Y$, respectively,
   in their kernel;
   we form tensor products of them, which yields characters with kernel $D$,
   and apply the LLL algorithm to them.

   This yields all missing irreducibles of $H$:
   The degree squares of all now known irreducibles sum up to the order
   of $H$, which means that no more class splitting occurs.
   
   Finally, we compute the power maps (and thus the element orders)
   of the character table of $H$.
   The result is a character table that is permutation equivalent to
   the character table that is stored in {\GAP}'s table library.
\end{itemize}

\section{Appendix: The character table of $3^{1+12}_+:6.\Suz.2$}%
\label{sect:norm3B}

\subsection{Overview}

The \texttt{3B} normalizer in $\M$ has the structure $3^{1+12}_+.2.\Suz.2$.
Its character table has been computed by Richard Barraclough and
Robert A.~Wilson, see \cite{BW07}.
In order to describe a reproducible construction of the character table
that does not assume the character table of $\M$,
we recompute this table.


Our approach is similar to that in \cite{BW07}:
The subgroup in question is a factor group of the split extension $H$ of
the extraspecial group $N = 3^{1+12}_+$ by $6.\Suz.2$,
and we compute the character table of this bigger group $H$.

%

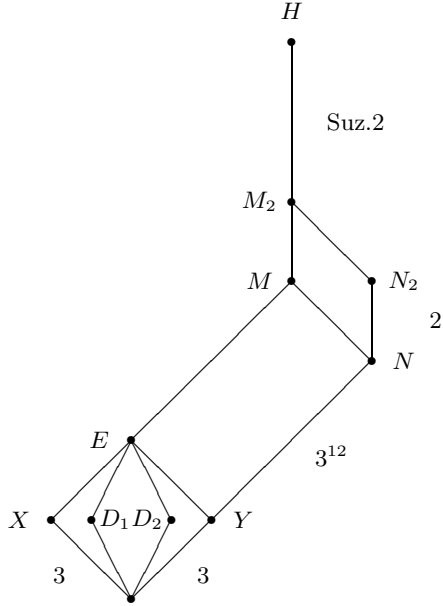
\begin{center}
\tthdump{\setlength{\unitlength}{3pt}}
\begin{picture}(50,80)
\put(20, 0){\circle*{1}}
\put(10,10){\circle*{1}} \put(6,10){\makebox(0,0){$X$}}
\put(15,10){\circle*{1}} \put(18,10){\makebox(0,0){$D_1$}}
\put(25,10){\circle*{1}} \put(22,10){\makebox(0,0){$D_2$}}
\put(30,10){\circle*{1}} \put(34,10){\makebox(0,0){$Y$}}
\put(20,20){\circle*{1}} \put(16,20){\makebox(0,0){$E$}}
\put(20, 0){\line(-1,1){10}}
\put(20, 0){\line(-1,2){5}}
\put(20, 0){\line(1,2){5}}
\put(20, 0){\line(1,1){10}}
\put(20,20){\line(-1,-1){10}}
\put(20,20){\line(-1,-2){5}}
\put(20,20){\line(1,-2){5}}
\put(20,20){\line(1,-1){10}}
\put(40,40){\circle*{1}} \put(36,40){\makebox(0,0){$M$}}
\put(40,50){\circle*{1}} \put(36,50){\makebox(0,0){$M_2$}}
\put(50,30){\circle*{1}} \put(54,30){\makebox(0,0){$N$}}
\put(50,40){\circle*{1}} \put(54,40){\makebox(0,0){$N_2$}}
\put(40,70){\circle*{1}} \put(40,74){\makebox(0,0){$H$}}
\put(20,20){\line(1,1){20}}
\put(30,10){\line(1,1){20}}
\put(50,30){\line(-1,1){10}} 
\put(50,40){\line(-1,1){10}} 
\put(40,40){\line(0,1){30}}  
\put(50,30){\line(0,1){10}}  
\put(11,3){\makebox(0,0){$3$}}
\put(29,3){\makebox(0,0){$3$}}
\put(45,18){\makebox(0,0){$3^{12}$}}
\put(48,60){\makebox(0,0){$\Suz.2$}}
\put(58,35){\makebox(0,0){$2$}}
\end{picture}
\end{center}

We have $H / N \cong 6.\Suz.2$.
Let $M$ be the normal subgroup of $H$ above $N$ such that
$H / M \cong 2.\Suz.2$ holds.
Then $M = X \times N$ for a normal subgroup $X$ of order $3$,
where the extension of $M / X$ by $H / M$ is split.
Let $N_2$ be the normal subgroup of $H$ above $N$ such that
$H / N_2 \cong 3.\Suz.2$ holds.
Then $M_2 = M N_2$ has the property $H / M_2 \cong \Suz.2$.

Let $Y = Z(N)$, of order $3$.
Then $E = X \times Y \cong 3^2$ is elementary abelian,
with diagonal normal subgroups $D_1$ and $D_2$.
It will turn out that the extensions of $M / D_1$ and $M / D_2$ by $H / M$
are non-split.

$H / Y$ is a subdirect product of $H / N$ and 
$H / E \cong 3^{12}.2.\Suz.2$,
we have $H / Y \cong (3 \times 3^{12}).2.\Suz.2$.

The group $H$ is a subdirect product of $H / N$ and $H / X$.

The \texttt{3B} normalizer in $\M$ is isomorphic to
one of the two factor groups $H / D_1$, $H / D_2$.
(The decision which of the two groups occurs as a subgroup of $\M$
appears in Section~\ref{sect:natcharM}.)

We use the following approach to compute the character table of $H$.

\begin{itemize}
\item
    Compute permutation generators \texttt{gensHmodX} of $H / X$,
    of degree $3^{13}$,
    see Section~\ref{HmodXpermgens}.
    The first two generators are standard generators of $2.\Suz.2$,
    the third generator lies in $M / X$.
\item
    Fetch permutation generators \texttt{gensHmodN} of $H / N_2$,
    of degree $5\,346$,
    and compute permutation generators \texttt{gensH} of $H$,
    of degree $1\,599\,669 = 1\,594\,323 + 5\,346$,
    see Section~\ref{Hpermgens}.
    The first two generators are standard generators of $6.\Suz.2$,
    the third generator lies in $N$.
\item
    Let {\MAGMA} compute the character table of $H$ from \texttt{gensH},
    see Section~\ref{computeHtable}.
\end{itemize}

\textbf{Remark:}
In an earlier construction (in September 2020),
we asked {\MAGMA} to compute the character table of $H / Y$,
and then used individual conjugacy tests for first setting up the class fusion
from $H$ to $H / Y$ and then determining the full character table head of $H$.
The computation of the missing irreducible characters of $H$ was then
not difficult, using character theoretic methods.
However, trying the same input files in 2024 showed many cases where
some of the required conjugacy tests did not finish in reasonable time.
Luckily, the automatic computation described in Section~\ref{computeHtable}
worked.

\subsection{A permutation representation of $H / X$}%
\label{HmodXpermgens}

The {\ATLAS} of Group Representations~\cite{AGRv3} contains
a faithful representation of $H$,
as a group of $38  \times 38$ matrices over the field with three elements.
The generating matrices are called \texttt{M3max7G0-f3r38B0.m1}, \ldots,
\texttt{M3max7G0-f3r38B0.m4}.
These matrices are block diagonal matrices with blocks of the lengths
$24$ and $14$.

\begin{verbatim}
    gap> info:= OneAtlasGeneratingSetInfo( "3^(1+12):6.Suz.2", Dimension, 38 );;
    gap> gens:= AtlasGenerators( info ).generators;;
    gap> Length( gens );
    4
    gap> ForAll( gens,
    >            m -> ForAll( [ 25 .. 38 ],
    >                         i -> ForAll( [ 1 .. 24 ],
    >                                      j -> IsZero( m[i,j] ) and
    >                                           IsZero( m[i,j] ) ) ) );
    true
\end{verbatim}

Here we will use just the lower right $14 \times 14$ blocks,
which generate the factor group $H / X \cong 3^{1+12}_+:2.\Suz.2$.
(Note that the $12$-dimensional irreducible representation of $2.\Suz.2$
over the field with $3$ elements respects a symplectic form,
and the embedding of $2.\Suz.2$ into $\Sp(12,3)$
--which is the automorphism group of the extraspecial group $3^{1+12}_+$--
yields a construction of the semidirect product $3^{1+12}_+:2.\Suz.2$
as a group of $14 \times 14$ matrices.)

Later we will construct a faithful permutation representation of $H$
as a subdirect product of $H / X$ and $H / N \cong 6.\Suz.2$.
Let \texttt{G} be the group generated by the $14 \times 14$ matrices.

When one deals with $H / X$, the fourth generator is redundant,
we will leave it out in the following.

\begin{verbatim}
    gap> mats:= List( gens, x -> x{ [ 25 .. 38 ] }{ [ 25 .. 38 ] } );;
    gap> List( mats, ConvertToMatrixRep );;
    gap> Comm( mats[3], mats[3]^mats[2] ) = Inverse( mats[4] );
    true
    gap> mats:= mats{ [ 1 .. 3 ] };;
\end{verbatim}

We use just the following facts.
The \texttt{G}-action on $GF(3)^{14}$ has six orbits,
from which we get a faithful permutation representation \texttt{homHtoHmodE}
of $H / E \cong 3^{12}.2.\Suz.2$ on $196\,560$ points
and a faithful representation \texttt{homHtoHmodX}
of a group of the structure $3^{1+12}_+.2.\Suz.2$ on $1\,594\,323$ points.
Since the image contains a subgroup $2.\Suz.2$,
the image is the split extension $H / X$ of $3^{1+12}_+$ by $2.\Suz.2$.


\begin{verbatim}
    gap> G:= GroupWithGenerators( mats );;
    gap> orbs:= ShallowCopy( OrbitsDomain( G, GF(3)^14 ) );;
    gap> Length( orbs );
    6
    gap> SortBy( orbs, Length );
    gap> List( orbs, Length );
    [ 1, 2, 196560, 1397760, 1594323, 1594323 ]
    gap> v:= [ 0, 0, 0, 0, 0, 0, 0, 0, 0, 0, 0, 0, 0, 2 ] * Z(3)^0;;
    gap> orb_large:= First( orbs, x -> v in x );;
    gap> Length( orb_large );
    1594323
    gap> Length( orb_large ) = 3^13;
    true
    gap> orb_large:= SortedList( orb_large );;
    gap> homHtoHmodX:= ActionHomomorphism( G, orb_large );;
    gap> represHmodX:= Image( homHtoHmodX );;
    gap> Size( represHmodX );
    2859230155080499200
    gap> Size( represHmodX ) = Size( CharacterTable( "2.Suz.2" ) ) * 3^13;
    true
\end{verbatim}

%
%
%
%

Next we show that the first two elements from the generating set
are standard generators of $2.\Suz.2$.
For that, we first show that these elements are preimages
of standard generators of $\Suz.2$,
by computing that the words in question lie in the centre of $2.\Suz.2$,
and then show that the elements satisfy the conditions of
standard generators of $2.\Suz.2$,
that is, the second generator has order $3$,
see the page on $\Suz$ in the {\ATLAS} of Group Representations~\cite{AGRv3}.


\begin{verbatim}
    gap> slp:= AtlasProgram( "Suz.2", "check" );;
    gap> prog:= StraightLineProgramFromStraightLineDecision( slp.program );;
    gap> res:= ResultOfStraightLineProgram( prog, gensHmodX );;
    gap> List( res, Order );
    [ 2, 1, 2, 1 ]
    gap> ForAll( gensHmodX{ [ 1, 2 ] }, x -> ForAll( res, y -> x*y = y*x ) );
    true
    gap> Order( gensHmodX[2] );
    3
\end{verbatim}

\subsection{A permutation representation of $H$}%
\label{Hpermgens}

The group $H$ is a subdirect product of $H / X$ and $H / N_2 \cong 3.\Suz.2$,
w.r.t. the common factor group $H / M_2 \cong \Suz.2$.
Since we know that our generators for $H / X$ are compatible with
standard generators of the factor group $2.\Suz.2$,
it is sufficient to form the diagonal product of our representation of
$H / X$ and a representation of $3.\Suz.2$ on standard generators.

\begin{verbatim}
    gap> 3suz2:= OneAtlasGeneratingSet( "3.Suz.2", NrMovedPoints, 5346 );;
    gap> 3suz2:= 3suz2.generators;;
    gap> omega:= [ 1 .. LargestMovedPoint( 3suz2 ) ];;
    gap> shifted:= omega + LargestMovedPoint( gensHmodX );;
    gap> pi:= MappingPermListList( omega, shifted );;
    gap> shiftedgens:= List( 3suz2, x -> x^pi );;
    gap> Append( shiftedgens, [ () ] );
    gap> gensH:= List( [ 1 .. 3 ], i -> gensHmodX[i] * shiftedgens[i] );;
    gap> NrMovedPoints( gensH );
    1599669
\end{verbatim}

The first two of the generators are standard generators of $6.\Suz.2$.
Note that we know already that they are preimages of standard generators
of $\Suz.2$,
it remains to show that they are elements $C$, $D$ where $D$ has order $3$
and $CDCDD$ has order $7$. 


\begin{verbatim}
    gap> Order( gensH[2] );
    3
    gap> Order( Product( gensH{ [ 1, 2, 1, 2, 2 ] } ) );
    7
\end{verbatim}

\subsection{Compute the character table of $H$}%
\label{computeHtable}

Now we let {\MAGMA} do the work.

\begin{verbatim}
    gap> H:= GroupWithGenerators( gensH );;
    gap> if CTblLib.IsMagmaAvailable() then
    >      mgmt:= CharacterTableComputedByMagma( H, "H_Magma" );
    >    else
    >      mgmt:= CharacterTable( "3^(1+12):6.Suz.2" );
    >    fi;
\end{verbatim}

This computation needed about three weeks of CPU time.


The result verifies the character table of $H$
that is available in the library.

\begin{verbatim}
    gap> IsRecord( TransformingPermutationsCharacterTables( mgmt,
    >        CharacterTable( "3^(1+12):6.Suz.2" ) ) );
    true
\end{verbatim}

This character table was used in Section~\ref{elements_3}.

\section{Appendix: The character table of $5^{1+6}_+.4.\J_2.2$}%
\label{sect:table_N5B}

The normalizer of a \texttt{5B} element in $\M$ has the structure
$5^{1+6}_+.4.\J_2.2$,
generators for this group as a permutation group of degree $78\,125$
are available in the {\ATLAS} of Group Representations~\cite{AGRv3}.
{\MAGMA}~\cite{Magma} can compute the character table from the group
within a few minutes,
the result turns out to be equivalent to the table that is available in
{\GAP}'s character table library.

\begin{verbatim}
    gap> g:= AtlasGroup( "5^(1+6):2.J2.4" );;
    gap> if CTblLib.IsMagmaAvailable() then
    >      mgmt:= CharacterTableComputedByMagma( g, "MN5B_Magma" );
    >    else
    >      mgmt:= CharacterTable( "5^(1+6):2.J2.4" );
    >    fi;
    gap> IsRecord( TransformingPermutationsCharacterTables( mgmt,
    >        CharacterTable( "5^(1+6):2.J2.4" ) ) );
    true
\end{verbatim}

This character table was used in Section~\ref{elements_5}.

\bibliographystyle{amsalpha}
\bibliography{manualbib.xml,../../../doc/manualbib.xml,../../atlasrep/doc/manualbib.xml}


\end{document}